\documentclass[11pt]{article}
\usepackage{amsfonts}
\usepackage{amssymb}
\newtheorem{lemma}{Lemma}[section]
\newtheorem{theorem}[lemma]{Theorem}
\newtheorem{proposition}[lemma]{Proposition}
\newtheorem{corollary}[lemma]{Corollary}
\def\dim{\mathop{\rm dim}\nolimits}
\def\deg{\mathop{\rm deg}\nolimits}

\def\char{\mathop{\rm char}\nolimits}
\def\rank{\mathop{\rm rank}\nolimits}

\def\ind{\mathop{\rm ind}\nolimits}

\def\Hilb{\mathop{\rm Hilb}\nolimits}
\def\Im{\mathop{\rm Im}\nolimits}

\def\cA{\mathop{\cal A{}}\nolimits}
\def\cW{\mathop{\cal W{}}\nolimits}
\def\cV{\mathop{\cal V{}}\nolimits}
\def\cP{\mathop{\cal P{}}\nolimits}
\def\cU{\mathop{\cal U{}}\nolimits}
\def\cO{\mathop{\cal O{}}\nolimits}

\def\cK{\mathop{\cal K{}}\nolimits}

\def\cZ{\mathop{\cal Z{}}\nolimits}
\def\cN{\mathop{\cal N{}}\nolimits}
\def\cM{\mathop{\cal M{}}\nolimits}
\def\cH{\mathop{\cal H{}}\nolimits}
\def\cX{\mathop{\cal X{}}\nolimits}
\def\cY{\mathop{\cal Y{}}\nolimits}
\def\cI{\mathop{\cal I{}}\nolimits}
\def\cR{\mathop{\cal R{}}\nolimits}

\def\cE{\mathop{\cal E{}}\nolimits}
\def\cF{\mathop{\cal F{}}\nolimits}

\def\cG{\mathop{\cal G{}}\nolimits}
\def\pf{\mbox{\bf Proof. }}

\newcommand{\N}{\mbox{$\mathbb N$}}
\newcommand{\C}{\mbox{$\mathbb C$}}
\newcommand{\R}{\mbox{$\mathbb P$}}
\newcommand{\Z}{\mbox{$\mathbb Z$}}
\newcommand{\Q}{\mbox{$\mathbb Q$}}
\newcommand{\n}{\mbox{${\displaystyle \frac{n}{2} }$}}
\newcommand{\m}{\mbox{${\displaystyle \Big[ \frac{n}{2} \Big] }$}}

\title{On the Irreducibility of Commuting Varieties\\ of Nilpotent
Matrices}
\author{Roberta Basili
\\ {\em Dipartimento di Matematica e Informatica,
Universit\`a di Perugia,} \\ {\em Via Vanvitelli 1, 06123 Perugia,
Italy, }
\small{\rm e-mail: basili@dipmat.unipg.it} }
\date{}
\begin{document}
\maketitle
\begin{abstract}
Given an $n\times n$ nilpotent matrix over an
algebraically closed field $K$, we prove some properties of the
set of all the $n\times n$ nilpotent matrices over $K$ which
commute with it. Then we give a proof of the irreducibility of the
variety of all the pairs $(A,B)$ of $n\times n$
nilpotent matrices over $K$ such that $[A,B]=0$ if
either $\char K=0$ or $\char K\geq \n $.
We get as a consequence a proof of the irreducibility of the local Hilbert
scheme of $n$ points of a smooth algebraic surface over $K$ if
either $\char K=0 $ or $\char K\geq \n $.
\end{abstract}
\section{Introduction}
Let $K$ be an algebraically closed field. T.S. Motzkin and O.
Taussky in \cite{Mo} and M. Gerstenhaber in \cite{Ge} proved that
for any  $n\in \N $ the variety of all the pairs $(A,B)$ of
$n\times n$ matrices over $K$ such that $[A,B]=0$ is irreducible.
R.W. Richardson in \cite{Ri} extended this result by showing that
if $\char K=0$ the variety of all the pairs of commuting elements of
a reductive Lie algebra over $K$ is
irreducible. D.I. Panyushev in \cite{Pa} studied this kind of
varieties in a more general context and, besides other results,
showed that if $\char K=0$ the variety of all the pairs of
commuting elements of a symmetric space of maximal rank is normal
(see also \cite{Bre}).
\newline Let $\cH (n,K)$ be the variety of all the pairs $(A,B)$ of $n\times n$ nilpotent
matrices over $K$ such that $[A,B]=0$. Recently V. Baranovsky 
proved in \cite{Bar} that $\cH (n,K)$ is irreducible if either $\char K=0$ or
$\char K> n$. The proof uses the
irreducibility of the local punctual Hilbert schemes of smooth
algebraic surfaces. In fact, let $\cX $ be an algebraic surface
over $K$ and let $\Hilb ^n \cX $ be the Hilbert scheme of $n$
points in $\cX $. Let $P$ be a nonsingular point of $\cX $ and let
$\Hilb ^n(\cO _P)$ be the fiber in $nP$ of the Hilbert-Chow
morphism from $\Hilb ^n \cX $ to $Sym^n(\cX )$. It parametrizes
the ideals of colenght $n$ of $\cO _P$. Let $\widehat {\cH }(n,K)$ be the
subvariety of $\cH (n,K)\times K ^n$ of all the triples $(A,B,v)$
such that $\dim \langle A^iB^jv\  :\  i,j=0,\ldots ,n-1\rangle
=n$. H. Nakajima in \cite{Na} showed that there exists a morphism
from $\widehat {\cH }(n,K)$ onto $\Hilb ^n (\cO _P)$ whose fibers are the
orbits of $\widehat {\cH } (n,K)$ with respect to the action of GL$(n,K)$. If
$x,y$ are local coordinates in $P$ and we regard $\cO _P$ as a
subset of $K[[x,y]]$, it associates to $(A,B,v)$ the ideal of all
$g\in \cO _P$ such that $g(A,B)v=0$. J. Brian\c{c}on in \cite{Bc}
proved that $\Hilb ^n(\cO _P) $ is irreducible if $\char K=0$ 
(see also \cite{Gra}). In
\cite{Ia} this was extended by A.A. Iarrobino to the case $\char
K> n$. \newline In this paper we give a more
elementary proof of the irreducibility of $\cH (n,K)$ if
either $\char K=0$ or $\char K\geq \n $. It consists of showing
that the subset of $\cH (n,K)$ of all $(A,B)$ such that $B$ is
regular, that is $\rank B=n-1$, is dense. For this purpose we fix
an $n\times n$ nilpotent matrix $B$  over $K$ and study the Jordan
canonical form of the $n\times n$ nilpotent matrices over $K$
which commute with $B$. If $A$ commutes with $B$ we determine
some properties of the Jordan canonical form of $B+f(A)$ for some
suitable $f\in K[x]$ such that $x\;  \mid \;  f$.
 We use the form of the centralizer of $B$ when $B$ is in
Jordan canonical form and the irreducibility of the variety of its
nilpotent elements. The main step of the proof of the irreducibility
of $\cH (n,K)$ is the following
result: if either $\char K=0$ or $\char K\geq \n $ and $m\in \N $ is such
that $\n +2< m\leq n$ the subvariety of $\cH (n,K)$ of all the
pairs $(A,B)$ such that $\rank A\geq n-2$ and  $\  \n +1<\ind A\leq m$
(where $\ind A$ is the index of nilpotency of $A$) is irreducible.
In the proof of this result we use the same map introduced by
Nakajima in \cite{Na} from an open subset of  $\widehat {\cH }
(n,K)$ to the grassmannian
of all the subspaces of $K[x,y]/(x,y)^n$
of codimension $n$. This result is not true for some $n,\ p\in \N $
such that $p< \n $ and $\char K=p$.\newline
As a consequence we get a proof of the
irreducibility of $\Hilb ^n(\cO _P) $ for
algebraically closed fields $K$ such that either $\char K=0$ or
$\char K\geq \n $.
\section{Preliminaries}
If $R$ is a ring and $n', n''\in \N$ let $\cM (n'\times n'',R)$, 
$\cM (n',R)$ and
$\cN (n',R)$ be the varieties of all the $n'\times n''$ matrices, of
all the $n'\times n'$ matrices and of all the $n'\times n'$ nilpotent
matrices respectively over $R$. Let $J_n$ be the nilpotent Jordan block of
order $n$ over $K$. In this section we will not use the hypothesis 
that $K$ is algebraically closed.
\begin{proposition}\label{1}
$\cN (n,K)$ is irreducible of dimension $n^2-n$.\end{proposition}
\pf If $\cN _1 (n,K)$ is the subspace of $\cN (n,K)$ of the
matrices whose entries of indices $(i,j)$ such that $j-i\neq 1$
are 0, the morphism from GL$(n,K) \times \cN _1 (n,K)$ to $\cN
(n,K)$ defined by $(G,A)\mapsto G^{-1}AG$ is surjective. Since the
centralizer of $J_n$ has dimension $n$, the dimension of the orbit
of $J_n$ is $n^2-n$,  which shows the claim. $\square$\newline
\newline
Let $B\in \cN (n,K)$ be fixed, let $\cN _B=\{ A\in \cN (n,K)\, :\,
[A,B]=0\} $ and let $u_1\geq \ldots \geq u_t$ be the orders of the
Jordan blocks of $B$. We choose a basis  $\Delta _B=\{
v_h^{u_h-1},\ldots ,v_h^0\  : \ h=1,\ldots ,t \} $  of $K^n$ with
respect to which $B$ is in Jordan canonical form. If $F$ is an extension field of $K$ 
and $ X\in \cM (n,F)$  we regard 
the matrix which represents $X$ with respect to $\Delta _B$ as a block 
matrix $(X_{hk})$, $h,k=1,\ldots ,t$,
where $X_{hk}\in \cM (u_h\times u_k,F)$. If $m\in \N$, $\ l,
l'=1,\ldots ,m$, and $(X_{ll'})$ is a block matrix, let $X_{(l)}$
and $X^{(l)}$ be respectively the $l$-th row and the $l$-th column
of blocks of $(X_{ll'})$ for $l=1,\ldots ,m$.
\begin{lemma} \label{2}
{\rm (H.W. Turnbull and A.C. Aitken, 1931)}. If $A\in \cM (n,K)$
we have $[A,B]=0$ if and only if for $1\leq k\leq h\leq t$ the
matrices $A_{hk}$ and $A_{kh}$ have the following form:
$$A_{hk}=\pmatrix{0 &\ldots & 0 & a_{hk}^1 & a_{hk}^2 & \ldots &
a_{hk}^{u_h} \cr \vdots &  &  & 0 & a_{hk}^1 & \ddots & \vdots \cr
\vdots &  &  &  & \ddots & \ddots & a_{hk}^2 \cr 0 & \ldots &
\ldots & \ldots & \ldots & 0 & a_{hk}^1 \cr},$$
$$A_{kh}=\pmatrix{a_{kh}^1 & a_{kh}^2 & \ldots & a_{kh}^{u_h} \cr
0 & a_{kh}^1 & \ddots & \vdots \cr \vdots & \ddots & \ddots &
a_{kh}^2\cr \vdots &  & 0 & a_{kh}^1 \cr \vdots &  &  & 0 \cr
\vdots &  &  & \vdots \cr 0 & \ldots & \ldots & 0 \cr}$$ where for
$u_h=u_k$ we omit the first $u_k-u_h$ columns and the last
$u_k-u_h$ rows respectively.
\end{lemma}
\pf See \cite{Tu}, \cite{Ga}, \cite{Bas}. $\square$\newline
\newline
Let $q_0=0$ and let $q_{\alpha }\in \{ 1,\ldots ,t\} $, $\alpha
=1,\ldots ,\hat t$, be such that $u_h=u_{h+1}$ if $q_{\alpha
-1}+1\leq h<q_{\alpha }$, $u_{q_{\alpha }}\neq u_{q_{\alpha }+1}$,
$q_{\hat t}=t$. If $(A_{hk})$, $h,k=1,\ldots ,t$, has the form of
Lemma \ref{2}, for  $\alpha, \beta = 1,\ldots ,\hat t $ let
$$\overline{A}_{\alpha \beta }=(a_{hk}^1),\qquad q_{\alpha
-1}+1\leq h\leq  q_{\alpha},\quad q_{\beta -1}+1\leq k\leq
q_{\beta }.$$
\begin{lemma} \label{3}
If $A\in \cM (n,K)$ is such that $[A,B]=0$ then $A\in \cN _B$ if
and only if $\overline{A}_{\alpha \alpha}$ is nilpotent for
$\alpha =1,\ldots ,\hat t$, hence $\cN _B$ is irreducible.
Moreover if $A\in \cN _B$ it is possible to choose $\Delta _B$ such
that $\overline A _{\alpha \alpha }$ is upper triangular for
$\alpha =1,\ldots ,\hat t$.
\end{lemma}
\pf For $l=0,\ldots ,u_1-1$ let $U^l=\langle v_h^l\ : \ h=1,\ldots
,t,\ u_h-1\geq l\rangle $; then $K^n={\displaystyle \bigoplus
_{l=0}^{u_1-1}} U^l$ and $A(U^l)\subseteq {\displaystyle \bigoplus
_{i=l}^{u_1-1}}U^i$. For $v\in K^n$ let $v={\displaystyle \sum
_{l=0}^{u_1-1}}v^{(l)}$ where $v^{(l)}\in U^l$ and let $L_{A,l}\,
:\, U^l\to U^l$ be defined by $L_{A,l}(v)=A(v)^{(l)}$. Then $A$ is
nilpotent if and only if $L_{A,l}$ is nilpotent for $l=0,\ldots
,u_1-1$. For $l=0,\ldots ,u_1-1$ let $\gamma _l\in \{ 1,\ldots
,\hat t\} $ be such that $l\leq u_{q_{\gamma _l}}-1$,
$l>u_{q_{\gamma _l+1}}-1$ if $\gamma _l\neq \hat t$. Then the
matrix of $L_{A,l}$ with respect to the basis $\{ v_h^l\  : \
h=1,\ldots ,t,\ u_h-1\geq l\} $ is the lower triangular block
matrix $(\overline A _{\alpha \beta})$, $\alpha ,\beta =1,\ldots
,\gamma _l$, which is nilpotent if and only if $\overline A _{\alpha
\alpha }$ is nilpotent for $\alpha =1,\ldots ,\gamma _l$. This proves
the first claim and, by Proposition \ref{1}, the second one.
\newline For $v\in K^n$ let $v={\displaystyle \sum _{\alpha = 1}^{\hat t}
v_{(\alpha )}}$ where $v_{(\alpha )}\in \langle v_h^l\  :\
h=q_{\alpha -1}+1,\ldots ,q_{\alpha },\  l=0,\ldots ,u_{q_{\alpha
}}-1\rangle $. For $\alpha =1,\ldots ,\hat t$ and $l=0,\ldots
,u_{q_{\alpha }}-1$ let $U^l_{\alpha}=\langle v_h^l\ :\
h=q_{\alpha -1}+1,\ldots ,q_{\alpha }\rangle $ and let $
L_{A,\alpha ,l}\ :\ U^l_{\alpha}\to U^l_{\alpha } $ be defined by
$L_{A,\alpha ,l}(v)=L_{A,l}(v)_{(\alpha )}$. Then the matrix of $L
_{A,\alpha ,l}$ with respect to the basis $\{ v_h^l\ : \
h=q_{\alpha -1}+1,\ldots ,q_{\alpha }\} $ is $\overline A _{\alpha
\alpha }$. Hence for $\alpha =1,\ldots ,\hat t$ there exists $\{
c(A,\alpha )_{hk}\in K\  :\  h,k=q_{\alpha -1}+1,\ldots ,q_{\alpha
}\} $ such that if we set
$$w_h^l=\sum_{k=q_{\alpha -1}+1}^{q_{\alpha
}}c(A,\alpha )_{hk}v_k^l$$ for $h=q_{\alpha -1}+1,\ldots
,q_{\alpha }$ we have that the matrix of $L_{A,\alpha ,l}$ with
respect to the basis $\{ w_h^l\  :\ h=q_{\alpha -1}+1,\ldots
,q_{\alpha }\} $ is upper triangular for $l=0,\ldots ,u_{q_{\alpha
}}-1$. Then the basis $\{ w_h^{u_h-1},\ldots ,w_h^0\ : \
h=1,\ldots ,t\} $ has the required property. $\square$
\newline  \newline
We denote by $\cN  _{B,\Delta _B}$ the subspace of $\cN _B$ of all
$A$ such that $\overline A_{\alpha \alpha }$ is upper triangular
for $\alpha =1,\ldots ,\hat t$.
\newline Let $n_i\in
\{1,\ldots ,t\} $, $i=1,\ldots ,r_B$, be such that $n_1=1$,
$u_{n_i}-u_{n_{i+1}}\geq 2$, $u_{n_i}-u_{n_{i+1}-1}\leq 1$ for
$i=1,\ldots ,r_B-1$, $u_{n_{r_B}}-u_t\leq 1$.
\begin{proposition}\label{4}
There exists a non-empty open subset of $\cN _B$ such that if $A$
belongs to it we have $\rank A = n-r_B$.
\end{proposition}
\pf By Lemma \ref{3} it is sufficient to prove that there exists a
non-empty open subset $\cZ $ of $\cN _{B,\Delta _B}$ such that
$\rank A=n-r_B$ for $A\in \cZ $. \newline If $r_B=1$ there exists
$A\in \cN _B$ regular such that $B$ is a power of $A$ and if $A\in
\cN _{B,\Delta _B}$ is regular then for $m=1,\ldots ,u_1-1$ we have
$$v_1^{u_1-j},\ldots ,v_t^{u_t-j}\in \langle
Av_1^{u_1-j-1},Av_2^{u_2-j},\ldots ,Av_t^{u_t-j}\  :\  j=1,\ldots
,m\rangle $$ for $j=1,\ldots ,m$. \newline Let $r_B>1$. For
$i=1,\ldots ,r_B$ and $A\in \cN _{B,\Delta _B}$ we define $V_{A,i}$
to be the following subset of $K^n$:
$$\langle
Av_{n_{i'}}^{u_{n_{i'}}-j-1},Av_{n_{i'}+1}^{u_{n_{i'}+1}-j},\ldots
,Av_{n_{i'+1}-1}^{u_{n_{i'+1}-1}-j}\  : \  j=1,\ldots ,i'-i+1,\
i'=1,\ldots ,i\rangle .$$ Then we have
$$V_{A,i}\subseteq \langle v_{n_{i'}}^{u_{n_{i'}}-j},
v_{n_{i'}+1}^{u_{n_{i'}+1}-j},\ldots
,v_{n_{i'+1}-1}^{u_{n_{i'+1}-1}-j}\  : \  j=1,\ldots ,i'-i+1,\
i'=1,\ldots ,i\rangle $$ for $i=1,\ldots ,r_B$. Moreover the open
subset $\cZ $ of $\cN _{B,\Delta _B}$ of all $A$ such that in this
relation equal holds for $i=1,\ldots ,r_B$ and $\rank A\geq n-r_B$
is non-empty. In fact, if $\check A\in \cN _{B,\Delta _B}$ is such
that $\check A_{hk}=0$ if there doesn't exist $i\in \{ 1,\ldots
,r_B\} $ such that $h,k\in \{ n_i,\ldots ,n_{i+1}-1\} $ and for
$h',k'\in \{ n_i,\ldots ,n_{i+1}-1\} $ the nilpotent matrix
$(\check A_{h'k'})$ is regular for $i=1,\ldots ,r_B$ then we have
$\check A\in \cZ $.\newline If $A\in \cZ $ we have
$$v_{n_i}^{u_{n_i}-1},v_{n_i+1}^{u_{n_i+1}-1},\ldots
,v_{n_{i+1}-1}^{u_{n_{i+1}-1}-1}\in V_{A,i}$$ for $i=1,\ldots
,r_B$, which implies that $\Im A$ is contained in the following
subspace of $K ^n$:
$$\langle Av_h^l\  :\ h=1,\ldots ,t,\ l=0,\ldots ,
u_h-1,\ (h,l)\neq (n_i,u_{n_i}-1)\ \mbox{\rm for } i=1,\ldots ,r_B
\rangle .$$ This shows that $\cZ $ has the required property.
$\square $
\section{On some automorphisms of $\cH (n,K)$} The aim of this
section is the proof of Proposition \ref{3.5} and Proposition
\ref{3.7}, from which we will get Corollary \ref{3.8}. For
this purpose, we first prove some properties of the ranks of the
matrices $(A^m)_{hk}$ for $A\in \cN _B$, $m\in \N$ and $h,k\in \{
1,\ldots ,t\} $.
\begin{lemma}\label{3.1} Let $A\in \cN _{B,\Delta
_B}$ and let $\{ i_1,\ldots ,i_r\} \subseteq \{ 1,\ldots ,t\} $ be
maximal such that $i_1<\ldots <i_r$, $u_{i_1}-u_{i_r}\leq 1$. If
$m\in \{ 2,\ldots ,r\} $ and $i,j\in \{ i_1,\ldots ,i_r\} $ we
have
$$\rank (A^m)_{ij}\leq  \left\{  \begin{array}{ll} u_j-1 & \mbox{\em if
}\  0\leq j-i\leq m-1\\ u_j-2 & \mbox{\em if }\   i-j\geq
r-m+1.\end{array}\right. $$
\end{lemma}
\pf If $j-i\geq 0$ and $m\in \{ 2,\ldots ,r\} $ we have $\rank
(A^m)_{ij}\leq u_j-1$ if there doesn't exist an $l\in \{ 1,\ldots
,t\} $ such that $(A^{m-1})_{il}$ and $A_{lj}$ have the first
column different from 0. If $i-j\geq 0$ and $m\in \{ 2,\ldots ,r\}
$ we have $\rank (A^m)_{ij}\leq  u_j-2$ if there doesn't exist an
$l\in \{ 1,\ldots ,t\} $ such that either $(A^{m-1})_{il}$ has the
first column different from 0 and $A_{lj}$ has the second column
different from 0 or the converse. Hence the claim is true for
$m=2$. Then by induction on $m$ and the previous observations we
get the claim. $\square $
\newline  \newline
Let $s_B$ be the maximum of the cardinalities of the subsets $\{
i_1,\ldots ,i_r\} $ of $\{ 1,\ldots ,t\} $ such that $i_1<\ldots
<i_r$, $\ u_{i_1}-u_{i_r}\leq 1$.\newline Since if $A\in \cN _B$
we have $(\overline A_{\alpha \alpha})^{s_B}=0$ for $\alpha
=1,\ldots ,\hat t$, by Lemma \ref{3.1} we get the following
corollary.
\begin{corollary} \label{3.2}
If $A\in \cN _B$, $\  i,j\in \{ 1,\ldots ,t\} $ and $|u_i-u_j|\leq 1$ we
have $\rank
(A^{s_B})_{ij}\leq \min \{ u_i,u_j\} -1$. $\square $
\end{corollary}
If $X$ is a matrix and $m\in \Z $ is such that $m\leq 0$
the claims $\rank X=m$ and $\rank X<m$ will mean $\rank
X=0$.\newline If $X_1\in \cM (p_1\times u,K)$ and $X_2\in \cM
(u\times p_2,K)$ have the following form:
$$X_i=\pmatrix{0 & \ldots & 0 & x_1^i & x_2^i & \ldots & x_{r_i}^i \cr
\vdots &  &  & 0 & x_1^i & \ddots & \vdots \cr \vdots &  &  &  &
\ddots & \ddots & x_2^i \cr \vdots &  &  &  &  & 0 & x_1^i \cr
\vdots &  &  &  &  &  & 0 \cr \vdots &  &  &  &  &  & \vdots \cr 0
& \ldots & \ldots & \ldots & \ldots & \ldots & 0\cr }$$ for
$i=1,2$, where $r_i\leq \min \{ p_i,u\} $, then $X_1X_2\in \cM
(p_1\times p_2,K)$ has the same form and $\rank (X_1X_2)=\rank
X_1+\rank X_2-u$.\newline Hence, by Lemma \ref{2}, for $m_1,m_2\in
\N $ and $A\in \cN _B$, since $A^{m_1+m_2}=A^{m_1}A^{m_2}$, we
have
\begin{equation}\label{e} \rank (A^{m_1+m_2})_{hk}\leq \max _{l\in \{ 1,\ldots ,t\}
}\{ \rank (A^{m_1})_{hl} +\rank (A^{m_2})_{lk} -u_l\}
.\end{equation}
Let $\widehat {\cN }_B$ be the subset of $\cN _B$ of
all $A$ such that if $i,j\in \{ 1,\ldots ,t\} $ and $|u_i-u_j|\leq
1$ we have $\rank A_{ij}\leq \min \{ u_i,u_j\} -1$. \newline Let
$i,j\in \{ 1,\ldots ,t\} $ be such that $i\leq j$. We associate to
the pair $(i,j)$ an $h(i,j)\in \N \cup \{ 0\} $ defined as
follows. If $u_i-u_j\leq 1$ we set $h(i,j)=0$. If $u_i-u_j\geq 2$
we define $h(i,j)$ to be the unique element of $\N$ such that
there exists $\{ k_l\  :\  l=0,\ldots ,h(i,j)\} \subseteq \{
1,\ldots ,t\} $ with the following properties: $i=k_0<k_1<\ldots
<k_{h(i,j)}$, $\  u_{k_l}-u_{k_{l+1}}\geq 2$ and
$u_{k_l}-u_{k_{l+1}-1}\leq 1$ for $l=0,\ldots ,h(i,j)-1$, $\ 0\leq
u_{k_{h(i,j)}}-u_j\leq 1$.
\begin{lemma} \label{3.3}
If $A\in \widehat {\cN }_B$, $i,j\in \{ 1,\ldots ,t\} $ and $m\in
\N $ we have
$$\rank (A^m)_{ij}\leq \min \{ u_i,u_j\} -m+h(\min \{ i,j\} ,
\max\{ i,j\} ).$$
\end{lemma}
\pf  The claim is true if $m=1$, hence we prove it by induction
on $m$. Let $m\geq 2$. \newline Let $j\leq i$. By \ref{e} we get
$$\rank (A^m)_{ij}\leq \max_{l\in \{ 1,\ldots ,t\} }\{ \rank
A_{il}+\rank (A^{m-1})_{lj}-u_l\} ,$$ hence it is sufficient to
prove that for $l=1,\ldots ,t$ we have
$$\rank A_{il}+\rank (A^{m-1})_{lj}-u_l\leq u_i-m+h(j,i).$$
Let $l\leq j$. Then we have $\rank A_{il}\leq u_i$ and equal may
hold only if $u_l-u_i\geq 2$. Moreover by the inductive hypothesis
we have $\rank (A^{m-1})_{lj}\leq u_j-m+1+h(l,j)$. Hence we get
$$\rank A_{il}+\rank (A^{m-1})_{lj}-u_l\leq
u_i+u_j-m+1+h(l,j)-u_l$$ where equal may hold only if $u_l-u_i\geq
2$. Since we have $u_j+h(l,j)\leq u_l$ and equal may hold only if
$u_j=u_l$ we get the claim.\newline Let $j\leq l\leq i$. Then as
above we have $\rank A_{il}\leq u_i$ and equal may hold only if
$u_l-u_i\geq 2$. By the inductive hypothesis we have $\rank
(A^{m-1})_{lj}\leq u_l-m+1+h(j,l)$, hence we get
$$\rank A_{il}+\rank (A^{m-1})_{lj}-u_l\leq u_i-m+1+h(j,l)$$
where equal may hold only if $u_l-u_i\geq 2$. Since $h(j,l)\leq
h(j,i)$ and equal may hold only if $u_l-u_i\leq 1$, we get the
claim.\newline Let $i\leq l$. Then we have $\rank A_{il}\leq u_l$
and equal may hold only if $u_i-u_l\geq 2$. By the inductive
hypothesis we have $\rank (A^{m-1})_{lj}\leq u_l-m+1+h(j,l)$,
hence we get
$$\rank A_{il}+\rank (A^{m-1})_{lj}-u_l\leq u_l-m+1+h(j,l)$$
where equal may hold only if $u_i-u_l\geq 2$. If $u_i=u_l$ then
$h(j,l)=h(j,i)$ and hence we get the claim. If $u_l+1=u_i$ then
$h(j,l)\leq h(j,i)+1$, hence we get the claim. Let $u_i-u_l\geq
2$. We have $h(j,l)\leq h(j,i)+h(i,l)+1$ and $u_l+h(i,l)+1\leq
u_i$, but not in both of these equal may hold. In fact
$u_l+h(i,l)+1=u_i$ only if $u_i-u_l=2$ and if $u_i-u_l=2$ we have
$h(j,l)\leq h(j,i)+h(i,l)$. Hence we get the claim.\newline If
$i\leq j$ we can repeat the same argument as above by using
$$\rank (A^m)_{ij}\leq \max_{l\in \{ 1,\ldots ,t\} }\{ \rank
(A^{m-1})_{il}+\rank A_{lj}-u_l\} .\  \  \square $$ 
\newline 
Let $m\in \N$ and let $\cN _m $ be the
subset of $\cM (n,K(x))$ of all $Y$ which commutes
with $B$ and such that
$$\rank Y_{ij}\leq \min \{ u_i,u_j\} -m+h(\min \{ i,j\} , \max \{
i,j\} )$$ for $i,j=1,\ldots ,t$. Let $q_0$ and $q_{\alpha }$, $\
\alpha =1,\ldots ,\hat t$, be as in Section 2 and let $t_m=\max \{
\alpha \in \{ 1,\ldots ,\hat t\} \  :\ u_{q_{\alpha }}-2\geq m\}
$.
\newline For $Y\in \cN _m$, $\ \alpha \in \{ 1,\ldots ,t_m\} $ and
$h,k\in \{ q_{\alpha -1}+1,\ldots ,q_{\alpha }\} $ let
$a_{hk}^m(Y)$ be the entry of $Y_{hk}$ of indices $(l,m+l)$, $\
l=1,\ldots , u_{q_{\alpha }}-m$, and let $d_{\alpha }^m(Y)=\det
(a_{hk}^m(Y))$. Let $\cV _m$ be the open subset of $\cN _m$ of all
$Y$ such that $d_{\alpha }^m(Y)\neq 0$ for $\alpha =1,\ldots
,t_m$. Let $\widetilde {\cV }_m$ be the subset of $\cV _m$ of all
$Y$ such that if $i,j\in \{ 1,\ldots ,t\} $ are such that
$|u_i-u_j|\geq 2$ then $Y_{ij}=0$. \newline For $Y\in \cN _m$ let
$\cR (Y)$ be the subset of $K$ of all $a$ such that all the
entries of $Y$ are regular in $a$ and let $\cE _m(Y)=\{ a\in K \
:\ d_{\alpha }^m(Y)(a)\neq 0\  \mbox{for }\ \alpha =1,\ldots
,t_m\} $.
\begin{lemma}\label{3.4}
If $m\in \N $ there exist an open subset $\cW _m$ of $\cV _m$ such
that  $\widetilde \cV _m \subset \cW _m$ and a rational map $\phi
_m\ :\ \cW _m \to \cN _m $ such that:\begin{itemize}
\item[1)] the first $m$ columns and the last $m$ rows of $\phi
_m(Y)_{hk}$ are 0 for $Y\in \cW _m$ and $h,k=1,\ldots ,t$;
\item[2)] if $Y\in \cW _m$ and $g$ is an entry
of $\phi _m(Y)$ then $$g\prod _{\alpha =1}^{t_m}d_{\alpha
}^m(Y)^{2[(\hat t-\alpha )+(\hat t-\alpha -1)+\cdots +(\hat
t-t_m)]}$$ is regular in $\cR (Y)$;
\item[3)] $\rank \phi _m (Y)(a)=\rank Y(a)$ for $Y\in \cW _m$ and
$a\in \cR (Y)\cap \cE _m (Y)$. \end{itemize}
\end{lemma}
\pf If $\alpha ,\beta \in \{ 1,\ldots ,\hat t\} $ are such that
$u_{q_{\alpha }}-u_{q_{\beta }}\geq 2$ and $Y\in \cV _m$, for
$i=q_{\alpha -1}+1,\ldots ,q_{\alpha }$ and $r=1,\ldots
,h(q_{\alpha },q_{\beta })$ let $C^i(Y,\beta ,r)$ be the submatrix
of $Y^{(i)}$ having as columns the columns of $Y^{(i)}$ of indices
$m+1,\ldots ,u_{q_{\beta }}+r$ and let $C_i(Y,\beta ,r)$ be the
submatrix of $Y_{(i)}$ having as rows the rows of $Y_{(i)}$ of
indices $u_{q_{\alpha }}-u_{q_{\beta }}-r+1,\ldots ,u_{q_{\alpha
}}-m$.  Let $r\in \{ 1,\ldots ,h(q_{\alpha },q_{\beta })\} $ and
let $j\in \{ q_{\beta -1}+1,\ldots ,q_{\beta }\} $. Let $Y\in \cV
_m$ be such that
\begin{equation}\label{e'}\rank Y_{ij},\ \rank Y_{ji}\leq
u_{q_{\beta }}-m+r\end{equation} \begin{equation}\label{e''}\rank
Y_{li},\  \rank Y_{il}\leq u_{q_{\alpha }}-m\end{equation} for
$i=q_{\alpha -1}+1,\ldots ,q_{\alpha }$ and $l=1,\ldots ,q_{\alpha
-1}$. Since $d_{\alpha }^m(Y)\neq 0$, it is possible to add to the
submatrix of $Y^{(j)}$ whose columns are the columns of $Y^{(j)}$
of indices $m-r+1,\ldots , u_{q_{\beta }}$ a suitable linear
combination over $K(x)$ of the matrices $C^i(Y,\beta ,r)$, $\
i=q_{\alpha -1}+1,\ldots ,q_{\alpha }$, getting a matrix $Y'$
which commutes with $B$ and such that $\rank (Y')_{ij}\leq
u_{q_{\beta }}-m+r-1$ for $i=q_{\alpha -1}+1,\ldots ,q_{\alpha }$.
The coefficients of this linear combination are a solution of a
Cramer system whose coefficients are entries of $Y$ and whose
matrix has $d_{\alpha }^m(Y)$ as determinant. Hence if $g$ is an
entry of $Y'$ we have that $gd_{\alpha }^m(Y)$ is regular in $\cR
(Y)$.
\newline We claim that $Y'\in \cN _m$. In fact, it is sufficient
to prove that for $i=q_{\alpha -1}+1,\ldots ,q_{\alpha }$ and
$l=q_{\alpha }+1,\ldots ,t$ we have
$$\rank C^i(Y,\beta ,r)_{(l)}\leq \min \{ u_l, u_{q_{\beta }}\} -m
+h(\min \{ l, q_{\beta }\} , \max \{ l, q_{\beta }\} ).$$ We have
$$\begin{array}{rcl}\rank C^i(Y,\beta ,r)_{(l)} & \leq & u_l-m+h(q_{\alpha }
,l)-u_{q_{\alpha }}+u_{q_{\beta }}+r\  \leq \\   & \leq &
u_l-m+h(q_{\alpha },l)-u_{q_{\alpha }}+u_{q_{\beta }}+h(q_{\alpha
},q_{\beta }).\end{array}$$ If $l\geq q_{\beta }$ then
$h(q_{\alpha },l)\leq h(q_{\alpha },q_{\beta })+h(q_{\beta },l)+1$
and $2h(q_{\alpha },q_{\beta })\leq u_{q_{\alpha }}-u_{q_{\beta
}}$, but not in both of these equal may hold. In fact, if
$2h(q_{\alpha },q_{\beta })= u_{q_{\alpha }}-u_{q_{\beta }}$ then
there exists $\{ k_{l'}\  :\  l'=0,\ldots ,h(q_{\alpha },q_{\beta
})\} $ such that $u_{k_{l'}}-u_{k_{l'+1}}=2$ for $l'=0,\ldots
,h(q_{\alpha },q_{\beta })-1$ and $k_0=q_{\alpha }$, $\
k_{h(q_{\alpha },q_{\beta })}=q_{\beta }$ and then $h(q_{\alpha
},l)=h(q_{\alpha },q_{\beta })+h(q_{\beta },l)$. If $l\leq
q_{\beta }$ we may use the same argument by reversing the role of
$q_{\beta }$ and $l$. This proves the claim.\newline If $Y\in
\widetilde {\cV }_m$  the coefficients of the above linear
combination of matrices are 0. Hence there exists an open subset
of $\cV _m$ which contains $\widetilde {\cV }_m$ such that if $Y$
belongs to it we get $Y'\in \cV _m$. We could repeat a similar
argument for the rows; in this case we would add to the submatrix
of $Y_{(j)}$ whose rows are the rows of $Y_{(j)}$ of indices
$1,\ldots ,u_{q_{\beta }}-m+r$ a suitable linear combination over
$K(x)$ of the matrices $C_i(Y,\beta ,r)$, $\ i=q_{\alpha
-1}+1,\ldots ,q_{\alpha }$. Let $O(\alpha ,\beta ,r)$ be the
operation on the rows and columns of an element $Y$ of $\, \cV _m$
which satisfies \ref{e'} and \ref{e''} for $i=q_{\alpha
-1}+1,\ldots ,q_{\alpha }$, $\ l=1,\ldots ,q_{\alpha -1}$ and
$j=q_{\beta -1}+1,\ldots ,q_{\beta }$ which consists on applying
the previous operation on the columns and then the corresponding
operation on the rows for $j=q_{\beta -1}+1,\ldots ,q_{\beta }$.
For $\alpha , \beta \in \{ 1,\ldots ,\hat t\} $ we define
$O(\alpha ,\beta , 0)$ to be the operation which doesn't change
the given matrix.
\newline If $\alpha ,\beta \in \{ 1,\ldots , \hat t\} $ and
$u_{q_{\alpha }}-u_{q_{\beta }}\geq 2$ then for $l=q_{\alpha }+1
,\ldots ,t$, $\ i=q_{\alpha -1}+1,\ldots , q_{\alpha }$ and
$r=1,\ldots ,h(q_{\alpha },q_{\beta })$ we also have
$$\rank C^i(Y,\beta ,r)_{(l)},\  \rank C_i(Y,\beta ,r)^{(l)}\leq
u_{q_{\beta }}-m+r-1$$ since $u_l+h(q_{\alpha }, l)+1\leq
u_{q_{\alpha }}$. Let $\cW _m$ be the open subset of $\cV _m$ of
all $Y$ such that it is possible to apply to $Y$ the following
operations in the order in which they are written:
$$O(\gamma ,\beta ,h(q_{\alpha },q_{\beta })),
\quad \gamma =\alpha ,\alpha -1,\ldots , 1, \   \beta =\alpha +1,\ldots
,\hat t, \  \alpha =1,\ldots ,t_m.$$ Let $\phi _m\ :\ \cW _m\to
\cN _m$ be the map which associates to $Y\in \cW _m$ the matrix
obtained by applying to $Y$ the previous operations in the given
order. Then $\phi _m$ has the required properties. $\square
$\newline \newline The previous lemma allows us to prove the main
results of this section.
\begin{proposition} \label{3.5}
If $A\in \cN _B$ and $m\in \N$ we have $\rank (A^{s_B})^m\leq
\rank B^m$ and $\rank (B+xA^{s_B})^m\leq \rank B^m$.
\end{proposition}
\pf The first claim is a consequence of the second one, hence let
us prove the second claim. It is sufficient to prove that there
exists a non-empty open subset of $\cN _B$ such that if $A$
belongs to it we have $\rank (B+xA ^{s_B})^m\leq \rank B^m$. By
Corollary \ref{3.2} and Lemma \ref{3.3} we get that $(B+xA
^{s_B})^m\in \cN _m$ for any $A\in \cN _B$. Moreover if $A=0$ we
have $(B+xA^{s_B})^m\in
\widetilde {\cV }_m$. Hence by {\em 1)} and {\em 3)} of Lemma
\ref{3.4} we get the claim. $\square$
\newline  \newline
Let $F$ be a field and for $p\in \N \cup \{ 0\} $ let $F[x]_p= \{ g\in
F[x]\ :\ \deg
g\leq p\} $.
\begin{lemma}\label{3.6}
\begin{itemize}
\item[i)] If $p_1,p_2\in \N \cup \{ 0\} $ the subset of $F[x]_{p_1}
\times F[x]_{p_2}$ of all $(g_1,g_2)$ such that $g_2\neq 0$ and
$g_1\not \in \sqrt{(g_2)}$ is open and if $p_2\in \N $ it is
non-empty.
\item[ii)] If $p_1,\ldots ,p_r\in \N \cup \{ 0\} $, $\  r>1$ and $F$ is
algebraically closed
the subset of $F[x]_{p_1}\times \cdots \times F[x]_{p_r}$ of all
$(g_1,\ldots ,g_r)$ such that there exists $i\in \{ 1,\ldots ,r\}
$ such that $\deg g_i=p_i$ and $1\in (g_1,\ldots ,g_r)$ is open
and non-empty.
\end{itemize}
\end{lemma}
\pf \begin{itemize}\item[{\em i)}] The claim is obvious if $p_1=0$
or $p_2=0$, hence we assume $p_1,p_2\in \N $. If $(g_1,g_2)\in
F[x]_{p_1}\times F[x]_{p_2}$ we have $g_1\in \sqrt{(g_2)}$ if and
only if $(g_1)^{p_2}\in (g_2)$. If $g_2\neq 0$ this happens if and
only if there exists $[a,h]\in \R (F\times F[x]_{p_1p_2})$ such
that $a(g_1)^{p_2}=hg_2$, which shows the claim.
\item[{\em ii)}] Let $x,y$ be homogeneous coordinates of $\R ^1
(F)$. For $p\in \N \cup \{ 0\} $ and $g=a_0+a_1x+\cdots +a_px^p\in
F[x]_p$ let $g_0=a_0y^p+a_1xy^{p-1}+\cdots +a_px^p\in F[x,y]$. If
$(g_1,\ldots ,g_r)\in F[x]_{p_1}\times \cdots \times F[x]_{p_r}$
and there exists $i\in \{ 1,\ldots ,r\} $ such that $\deg g_i=p_i$
then $1\in (g_1,\ldots ,g_r)$ if and only if there doesn't exist
$[\overline x, \overline y]\in \R ^1(F)$ such that
$(g_j)_0(\overline x,\overline y)=0$ for $j=1,\ldots ,r$, which
shows the claim. $\square $
\end{itemize}
\begin{proposition}\label{3.7}
Let $\{ i_1,\ldots ,i_{s_B}\} \subseteq \{ 1,\ldots ,t\} $ be such
that $i_1<\ldots <i_{s_B}$, $\  u_{i_1}-u_{i_{s_B}}\leq 1$ and
$u_{i_1}\geq 2$. If $A\in \cN _B$ is such that at least one of the
determinants of the minors of $(B+xA ^{s_B})^{u_{i_1}-1}$ of order
$\rank B^{u_{i_1}-1}$ has the maximum possible degree there
exists $a(A)\in K$ such that
\begin{itemize}\item[a)] $\rank
(B+a(A)A^{s_B})^{u_{i_1}-1}<\rank B^{u_{i_1}-1}$.
\end{itemize}
Moreover there exists a non-empty open subset $\cU _{i_1}$ of $\cN
_B$ such that if $A\in \cU _{i_1}$ there exists $a(A)\in K$ such
that a) and the following conditions hold:
\begin{itemize}\item[b)]
$\rank (B+a(A)A^{s_B})\geq \rank B-1$;
\item[c)]$\rank (B+a(A)A^{s_B})^l=\rank B^l$
for $l=u_{i_1},\ldots ,u_1$.
\end{itemize}
\end{proposition}
\pf By Corollary \ref{3.2} and Lemma \ref{3.3} we have
$(B+xA^{s_B})^{u_{i_1}-1}\in \cN _{u_{i_1}-1}$ for any $A\in \cN
_B$. Let $\cW _{u_{i_1}-1}$ and $\phi _{u_{i_1}-1}$ be as in Lemma
\ref{3.4} and let $\overline {\cW }_{i_1}$ be the open subset of
$\cN _B$ of all $A$ such that $(B+xA^{s_B})^{u_{i_1}-1}\in \cW
_{u_{i_1}-1}$. For $A\in \overline {\cW }_{i_1}$ let $d_{i_1}(A)$
be the determinant of the submatrix of $(\phi
_{u_{i_1}-1}((B+xA^{s_B})^{u_{i_1}-1})_{hk})$, $\  h,k=1,\ldots
,t$,  obtained by taking the first $u_h-(u_{i_1}-1)$ rows and the
last $u_k-(u_{i_1}-1)$ columns of $\phi
_{u_{i_1}-1}((B+xA^{s_B})^{u_{i_1}-1})_{hk}$ for $h,k=1,\ldots
,\max \{ i\in \{ 1,\ldots,t\} \  :\  u_i\geq u_{i_1}\} $. Let
${\displaystyle m_{i_1}=\sum _{i=1}^t\max \{ u_i-(u_{i_1}-1),0\}
}$ and let
$$d_{i_1}'(A) = d_{i_1}(A) {\displaystyle \prod _{\alpha
=1}^{t_{u_{i_1}-1}}}d_{\alpha
}^{u_{i_1}-1}((B+xA^{s_B})^{u_{i_1}-1})^{2m_{i_1}[(\hat t-\alpha
)+(\hat t-\alpha -1)+\cdots +(\hat t-t_{u_{i_1}-1})]}.$$ By {\em
2)} of Lemma {\ref{3.4} we have $d_{i_1}'(A)\in K[x]$. We claim
that the subset $\overline {\cU }_{i_1}$ of $\overline {\cW
}_{i_1}$ of all $A$ such that $d_{i_1}'(A)\neq 0$ and
$$\prod _{\alpha =1}^{t_{u_{i_1}-1}}d_{\alpha
}^{u_{i_1}-1}((B+xA^{s_B})^{u_{i_1}-1})\not \in
\sqrt{(d_{i_1}'(A))}$$ is non-empty. In fact, let $\check A\in \cN
_B$ be such that $\check A_{hk}=0$ if $h\not \in \{ i_1,\ldots ,
i_{s_B}\} $ or $k\not \in \{ i_1,\ldots , i_{s_B}\} $ and for $h',
k'=i_1,\ldots ,i_{s_B}$ we have $(B_{h'k'})=\check a_1(\check
A_{h'k'})^{s_B}+\check a_2(\check A_{h'k'})^{s_B+1}$ where
$\check a_1,\check a_2\in
K\setminus \{ 0\} $. We have $d_{\alpha } ^{u_{i_1}-1}((B+x\check
A^{s_B})^{u_{i_1}-1})=1$ for $\alpha =1,\ldots ,t_{u_{i_1}-1}$.
Since $(B+x\check A^{s_B})^{u_{i_1}-1}\in \widetilde {\cV
}_{u_{i_1}-1}$ we have $\check A\in \overline {\cW }_{i_1}$. If we
set $a(\check A)=-\check a_1$ then $\check A$ satisfies {\em a)} and
$-\check a_1$ is the only element of $K$ with this property. Hence by
{\em 3)} of Lemma \ref{3.4} we have $\deg d_{i_1}'(\check A) \geq
1$ and then $\check A\in \overline {\cU }_{i_1}$.\newline By {\em
i)} of Lemma \ref{3.6} the subset $\overline {\cU }_{i_1}$ is open
in $\cN _B$. If $A\in \overline {\cU }_{i_1}$ then all the
determinants of the minors of $(B+xA^{s_B})^{u_{i_1}-1}$ of order
$\rank B^{u_{i_1}-1}$ have a common root, hence by {\em ii)} of
Lemma \ref{3.6} we get that if $A\in \cN _B$ is such that at least
one of these determinants has the maximum possible degree then
they have a common root. This proves the first claim. \newline If
we set $a(\check A)=-\check a_1$ then $\check A$ satisfies also {\em b)}
and {\em c)}. Hence the subset $\cU _{i_1}$ of $\overline {\cU }_{i_1}$
of all $A$ such that if
$\rank B>1$ there exists a minor of $B+xA^{s_B}$ of order $\rank
B-1$ and determinant $d(A)$ and, if $i_1>1$, there exists a minor
of $(B+xA^{s_B})^l$ of order $\rank B^l$ and determinant $d^l(A)$
for $l=u_{i_1},\ldots ,u_1-1$ such that if $\rank B>1$ we have
$$d(A)\prod _{\alpha =1}^{t_{u_{i_1}-1}}d_{\alpha
}^{u_{i_1}-1}((B+xA^{s_B})^{u_{i_1}-1})\not \in
\sqrt{(d_{i_1}'(A))}$$
and if $i_1>1$ we have
$$d(A) \prod _{l=u_{i_1}}^{u_1-1}d^l(A)\prod _{\alpha
=1}^{t_{u_{i_1}-1}}d_{\alpha
}^{u_{i_1}-1}((B+xA^{s_B})^{u_{i_1}-1})\not \in
\sqrt{(d_{i_1}'(A))}$$
is non-empty. Hence by {\em i)} of Lemma \ref{3.6} and Proposition
\ref{3.5} we get the second claim.
$\square $
\begin{corollary}\label{3.8}
Let $\cA $ be an open subset of $\cH (n,K)$, let $(A,B)\in \cA $ 
and let $\{ i_1,\ldots ,i_{s_B}\} \subseteq \{
1,\ldots ,t\} $ be such that $u_{i_1}-u_{i_{s_B}}\leq 1$. There exist 
$(A', B') \in \cA $ and $f\in
K[x]$ such that $x\;  \mid \;  f$, $B'+f(A')$ has $s_{B'+f(A')}$ Jordan blocks
of order 1 and if we denote by $u'_1,\ldots ,u'_{t'}$ the orders of the
Jordan blocks of $B'+f(A')$ we have
 $t'\geq t$ and $u'_i=u_i$ for $i=1,\ldots ,i_1-1$.
\end{corollary}
\pf Let $\cU _{i_1}$ be as in Proposition \ref{3.7}, let
$A_1\in \cU _{i_1}$ be such that $(A_1,B)\in \cA $ and let
$a(A_1)$ be as in Proposition \ref{3.7}. We set $a_1=a(A_1)$ and
$s_1=s_B$. We now set $B$ to be the matrix $B+a_1(A_1)^{s_1}$. If
it is possible, let us choose a subset $\{ i_1 ,\ldots , i_{s_B}\}
$ of $\{1,\ldots ,t\} $ maximal such that $u_{i_1}-u_{i_{s_B}}\leq
1$ and $u_{i_1}>1$ and let $\cU _{i_1}$ be as in Proposition
\ref{3.7}. Let $A_2\in \cU _{i_1}$ be such that
$(A_2,B-a_1(A_2)^{s_1})\in \cA $ and let $a(A_2)$ be as in
Proposition \ref{3.7}. We set $a_2=a(A_2)$ and $s_2=s_B$. We now
set $B$ to be the matrix $B+a_2(A_2)^{s_2}$. If it is possible,
let us repeat the previous argument, finding a suitable matrix
$A_3$ such that $(A_3,B-a_1(A_3)^{s_1}-a_2(A_3)^{s_2})\in \cA $.
If we repeat this argument, by Proposition \ref{4} at the end we
get a matrix $B'\in \cN (n,K)$ and a matrix $A'\in \cN _{B'}$
which have the required properties. $\square $
\begin{proposition}\label{3.8'}
Let $u_1-u_2=2$ and $u_i-u_{i+1}\geq 2$ for $i=2,\ldots ,t-1$.
There exists a non-empty open subset of $\cN _B$ such that if $A$
belongs to it there exists $a(A)\in K$ such that $\ind
(B+a(A)A)<u_1$, $\rank (B+a(A)A)=n-t$.
\end{proposition}
\pf 
If $A\in \cN _B$ by Lemma \ref{3} and Lemma \ref{3.3} 
the only entry of the matrix which represents $(B+xA )^{u_1-1}$ 
with respect to $\Delta _B$ which may be different from 0 is the entry 
of indices $(1,u_1)$; let us denote it by $h(A)$.\newline 
By {\em i)} of Lemma \ref{3.6} it is sufficient to prove that the 
subset $\cK $ of $\cN _B$ of all $A$ such that $h(A)\neq 0$
and there exists a minor $M$ of $B+xA $ of order $n-t$ such that 
$\det M\not \in \sqrt{(h(A))}$ is non-empty.\newline 
If $X\in \cM (n,K(x,y))$ let $(X_{12})'$ be the
submatrix of $X_{12}$ obtained by taking the first $u_2$ rows
and let $(X_{21})'$
be the submatrix of  $X_{21}$ obtained
by taking  the last $u_2$
columns.\newline 
Let $\check A \in \cN (n,K(y))$ be such that $(\check A
_{12})',(\check A _{21})'=I_{u_2}$, $\check A _{22}=yJ_{u_2}$, $\check
A_{ij}=0$ if $i,j\in \{ 1,\ldots ,t\} $ and $(i,j)\neq 
(1,2),(2,1),(2,2)$;
then $[\check A,B]=0$.\newline
Let us consider the following matrix of $\cM (2,K(x,y))$:
$$H=\pmatrix{x&1\cr 1& x+y\cr } $$
and for $k=1,\ldots ,u_1-1$ let $H^k=(h^k_{ij})$, $i,j=1,2$. 
Then we have
$((\check A+xB)
^k)_{ii}=h^k_{ii}
(J_{u_i})^k$, $\  (((\check A+xB)^k)_{ij})'=h^k_{ij}
(J_{u_2})^{k-1}$ for $i=1,2$, $\  j\in \{ 1,2\} \setminus \{ i\} $ and
$k=1,\ldots ,u_1-1$.
For $k=1,\ldots ,u_1-1$ and $i,j=1,2$ let us regard 
$h^k_{ij}(x,0)$ and $h^k_{ij}(0,y)$ as polynomials over $K(x)$ 
and $K(y)$ respectively. 
By induction
on $k$, we get that $h^k_{11}(x,0)$ is monic of degree $k$,
$\  h^k_{11}(0,y)$ is monic of degree $k-2$,
$\  h^k_{12}(0,y)$ and  $h^k_{21}(0,y)$ are monic of
degree $ k-1 $,  $\  h^k_{22}(0,y)$ is monic
of degree $k$ for $k=2,\ldots ,u_1-1$.\newline
If $a,b\in K$ we have $\rank (\check A(b)+aB)<n-t$ if and only if 
either $t>2$ and $a=0$ or $\det H(a,b)=0$. \newline
We claim that there exists a non-empty open subset $K'$ of $K$ 
such that if $b$ belongs to it we have 
$$(\det H(x,b),h^{u_1-1}_{11}(x,b))=1.$$ 
Since $\det H=x^2+yx-1$ if $\char K\neq 2$ this follows from the 
formula of the roots of $\det H $ as polynomial over $K(y)$. Let 
$\char K=2$. By {\em ii)} of Lemma \ref{3.6} it is sufficient to 
prove that 
$$(\det H(x,1),h^{u_1-1}_{11}(x,1))=1.$$
We have $\det H(x,1)=x^2+x+1$, hence the roots of $\det H(x,1)$ 
are the 3-th roots of 1 different from 1. If, for $k=1,\ldots 
,u_1-1$, we denote by $\widetilde H_k$ the matrix over $K(x)$ 
obtained by substituting in $H(x,1)^k$ the exponent of any power 
of $x$ with the remainder of the division of it by 3, we have
$$\widetilde H_2=\pmatrix{x^2+1 & 1 \cr 1 & x^2\cr },\quad
\widetilde H_3=\pmatrix{x & x^2+x \cr x^2+x & x^2\cr },\quad \widetilde  
H_4=\widetilde H_1.$$
Hence $h^k_{11}(x,1)$, for $k=1,\ldots ,u_1-1$, hasn't a root 
which is a 3-th root of 1 different from 1, which proves the 
previous claim if $\char K=2$. \newline
Let $K''$ be the non-empty open subset of $K$ of all $b$ such that  
$h^{u_1-1}_{11}(0,b)\neq 0$. Let $\overline b\in K'\cap K''$ and let 
$\overline a\in K\setminus \{ 0\} $ be such that $h^{u_1-1}_{11}(\overline
a,\overline
b)=0$.  Then $\ind (B+\overline a^{-1}\check A (\overline b))<u_1$, $\  \rank 
(B+\overline a^{-1}\check A (\overline b))=n-t$. Since
$h^{u_1-1}_{11}(x,\overline b)\neq 0$ 
we have
$h(\check A (\overline b))\neq 0$. Hence $\check A(\overline b)\in \cK $.
$\square $
\section{Proof of the irreducibility of $\cH (n,K)$}
For $q\in \Q
$ let $[q]=\min \{ q'\in \Z \  :\ q'\geq q\} $.
Let $$\widetilde {\cH }(n,K)=\{ (A,B)\in \cH (n,K)\  :\  \rank A\geq n-2\} . $$
For $m\in \Big\{ 3,\ldots , \m \Big\} $ let
$$\cH _m(n,K)=\{ (A,B)\in \cH (n,K)\ :\ \ind A\leq n-m+3\} $$ and let
$$\widetilde {\cH }_m(n,K)=\Big \{ (A,B)\in \widetilde {\cH }(n,K)\  :\  \n +1<\ind
A\leq n-m+3 \Big\} .$$
We first prove the irreducibility of the variety $\widetilde {\cH
}_m(n,K)$ for $m=3,\ldots ,\m $ under some conditions on the
characteristic of $K$.\newline
For $l\in \Big\{ 3,\ldots , \m
\Big\} $ let
$$J_{l,n}=\pmatrix{J_l & 0 \cr 0 & 0 \cr }.$$
\begin{lemma}\label{4.1}
Let $\cA $ be an open subset of $\widetilde {\cH }(n,K)$, let $(A,B)\in \cA $ 
and let $\n +1<\ind A<n$. There exist $(A',B')\in \cA $ and $f\in K[x]$ such that 
$x\; \mid \; f$ and $B'+f(A')$ is conjugate to $J_{n-\ind A +2,n}$. 
\end{lemma}
\pf 
Since $\{ A\} \times \cN _A$ is irreducible, by Proposition \ref{3.5} we may 
choose $B$ such that $A$ and $B$ are conjugate. We denote by 
$u_1$ and $u_2$, where $u_1-u_2>2$, the orders of their Jordan 
blocks.\newline
Let $\Delta _A$ be a basis of $K^n$ such that $A$ with respect to it is in Jordan
canonical form and for $X\in \cN _A$ let $(X_{ij}^{(A)})$, $i,j=1,2$ be the 
block matrix
which represents $X$ with respect to $\Delta _A$, where $X_{ij}^{(A)}$ is an 
$u_i\times u_j$
matrix for $i=1,2$.
By Lemma \ref{2} there exists $\widetilde f\in K[x]$ such that $x\; \mid \; \widetilde f$ and 
$(B+\widetilde f(A))_{11}^{(A)}$
is 0. Then $\rank (B+\widetilde f(A))\leq n-3$.
If we choose $B$ in a suitable open subset of $\cN _A$ we
have 
$$\rank ((B+\widetilde f(A))^h)_{12}^{(A)}, \  \rank
((B+\widetilde f(A))^h)_{21}^{(A)}=u_2+1-h,$$
$$\rank ((B+\widetilde f(A))^h)_{22}^{(A)}=u_2-h$$
and
$$\rank
((B+\widetilde f(A))^h)_{11}^{(A)}=u_2+2-h$$ 
for $h=2,\ldots , u_2+1$.
Hence
the difference between the orders of the first two Jordan blocks
of $B+\widetilde f(A)$ is greater than 1.
If we apply Corollary \ref{3.8} to $(A,B+\widetilde f(A))$ and the open subset 
$\{ (A'',B''+\widetilde f(A''))\ :\ 
(A'',B'')\in \cA \} $ by Proposition
\ref{4} we get that
there exists $(A',B')\in \cA $ and $f\in K[x]$ such that $x\; \mid \; f$ and, if we denote
by $u'_1,\ldots ,u'_{t'}$ the orders of the Jordan blocks of $B'+f(A')$, we have
$u'_1=u_2+2$ and $u'_i=1$ for $i=2,\ldots ,t'$.
This proves the claim. $\square$
\newline  \newline
For $l\in \Big\{ 3,\ldots ,\m \Big\} $ and $c\in K\setminus \{ 0\} $ let 
$$A_{l,n}(c)=\pmatrix{J_l & E_{l,n}(c)\cr E'_{l,n} & J_{n-l}}$$
where
$$E_{l,n}(c)=\pmatrix{c & 0 & \ldots & 0\cr 0 & \ldots & \ldots  &
0 \cr \vdots &  &  & \vdots \cr 0 & \ldots & \ldots & 0\cr },
\quad 
E'_{l,n}=\pmatrix{0 & \ldots & \ldots & 0\cr \vdots &  &  & \vdots \cr
 0 & \ldots & \ldots & 0\cr 0 & \ldots & 0 & 1 \cr}.$$
For $l\in \Big\{ 2,\ldots ,\m \Big\}  $ let
$$Z_l=\left\{ \begin{array}{ll} \{ (n-1,0),\ldots ,(0,0)\}  & \mbox{\rm if }
\  l=2\\  &  \\ 
\begin{array}{l} \{ (n-l+1,0), \ldots ,(l-2,0),
(l-3,1),\\ (l-3,0),(l-4,1),(l-4,0),\ldots ,(0,1),(0,0)\} \end{array} & 
\mbox{\rm if } \  l>2,\end{array}\right.$$
which has cardinality $n$.
We denote by $\{ e_1,\ldots ,e_n\} $ the canonical basis of $K^n$. 
\begin{lemma}\label{4.2} For $l\in
\Big\{ 3,\ldots , \m \Big\} $
the subset of $\cN _{J_{l,n}}$ of all $A$ such that 
$\rank A =n-2$, $\  \ind A=n-l+2$ is open and non-empty and if $A$ belongs to it then:
\begin{itemize}
\item[{\em a)}] there exist $c\in K\setminus \{ 0\} $, $\ \widehat f\in K[x]$ 
such that $\deg \widehat f\leq l-2$, $\  x\;  \mid \;  \widehat f$, $\  x^2\;  \nmid \;  \widehat f$ 
and a basis of $K^n$ such that the matrices which represent $A$ 
and $J_{l,n}$ with respect to it are $A_{l,n}(c)$ and $\widehat f(J_{l,n})$ 
respectively;
\item[{\em b)}] 
if $v\in K^n$ is such that $(J_{l,n})^{l-1}v\neq 0$ then 
$$\{ v,Av,\ldots ,A^{n-l+1}v,J_{l,n}v,\ldots ,(J_{l,n})^{l-2}v\} 
,\  \{ A^i (J_{l,n})^jv\ :\ (i,j)\in Z_l\} $$
are linearly independent.
\end{itemize}
\end{lemma}
\pf By Lemma \ref{2} and  Lemma \ref{3} if $A\in \cN 
_{J_{l,n}}$ there exists $g_A\in K[x]$ such that $x\; \mid \;  g_A$ 
and the submatrix of $A$ obtained by taking the first $l$ rows and 
columns is $g_A(J_l)$. For $A\in \cN _{J_{l,n}}$ let $A_0=A-g_A(J_{l,n})$ and let 
$A _0'$ be the submatrix of $A_0$ obtained by taking the rows and the columns 
of indices 
$1,l,\ldots ,n$. Since the rows and the columns of $A_0$ of indices $2,\ldots ,l-1$ are 0
we have $\ind A_0=n-l+2$ if and only if $ A _0'$, which is nilpotent, is regular.
By Lemma \ref{2} and Lemma \ref{3} the subset of $\cN _{J_{l,n}}$ of all $A$ such that 
$A _0'$ is regular is non-empty. Since $ A_0J_{l,n}=0$ we have $\ind A=n-l+2$ if and only if 
$\ind A_0=n-l+2$. If $A_0'$ is regular we have $\rank A=n-2$ if and only if $x^2\; \nmid \; g_A$.
Hence the subset of $\cN _{J_{l,n}}$ of all $A$ such that 
$\rank A=n-2$ and $\ind A=n-l+2$ is open and non-empty. Let $A$ 
belong to this subset. Then
$\{ e_1,\ldots ,e_l,  (A_0) ^{n-l}e_l,\ldots ,A_0 e_l\} $
is linearly independent and there exists $a_A\in K\setminus \{ 0\} $ such
that $ (A_0) ^{n-l+1}
e_l=a_Ae_1$. 
Then $\{ e_l,Ae_l,\ldots ,A^{n-l+1}e_l,J_{l,n}e_l,\ldots 
,(J_{l,n})^{l-2}e_l\} $ is linearly independent.
Since $AJ_{l,n}=g_A(J_{l,n})J_{l,n}$ and $ x^2\;  \nmid \; g_A$ we get b).\newline
Let $\gamma _1,\ldots ,\gamma _{l-1}\in K$, $\  \gamma _1\neq 0$,  be such
that $J_{l,n}=\gamma _1
g_A(J_{l,n}) +\cdots + \gamma _{l-1} g_A(J_{l,n})^{l-1}$ and let 
$\widehat f=\gamma _1x+\cdots +\gamma _{l-2}x^{l-2}$.  
Let $e_{l-i}'=g_A(J_{l,n})^ie_l$ for
$i=1,\ldots ,l-1$. If $c'\neq 0$ is the coefficient of $x$ in $g_A$ we have 
$e_1'=(c')^{l-1}e_1$ and hence if $c=a_A(c')^{1-l}$ we have $(A_0)^{n-l+1}e_l=ce_1'$. 
Then $\widehat f,\ c$ and the 
basis
$$\{ e_1',\ldots ,e_{l-2}',e_{l-1}'+\gamma _1^{-1}\gamma
_{l-1}e_1',e_l+\gamma _1 ^{-1}\gamma _{l-1}c^{-1}( A_0)
^{n-l}e_{l},( A_0)
^{n-l}e_l,\ldots , A_0e_l\} $$ 
have the property required in a).
$\square$
\newline  \newline
For $l\in \N \cup \{ 0\} $ and $l'\in \N $ let
$$F(l,l')=\left\{ \begin{array}{ll}y_1y+y_2y^2+\cdots +y_{l'}y^{l'} & 
\mbox{\rm if }
l=0\\
x_1x+x_2x^2 +\cdots +x_{l}x^{l}+y_1y+y_2y^2
+\cdots +y_{l'}y ^{l'} & \mbox{\rm if } l\neq 0,\end{array}\right.$$ 
which we consider as a polynomial of $\Z [y_1,\ldots ,y_{l'}][x,y]$ if $l=0$ and of $\Z
[x_1,\ldots ,x_{l},y_1,\ldots ,y_{l'}][x,y]$ if $l\neq 0$.\newline
For $k\in \N \cup \{ 0\} $ let $F(l,l',k)$ be
the polynomial obtained by substituting $x^iy^j$ with $y^{i+j}$ 
in $F(l,l')^k$ for $i,j\in \N $. Let $f_0(l,l',k)$ be the term of 
degree 0 of $F(l,l',k)$ and for $h\in \N $ let $f_h(l,l',k)$ and $f'_h 
(l,l',k)$ be the coefficients of $x^h$ and $y^h$ respectively in 
$F(l,l',k)$.
\begin{lemma}\label{A}
Let $l\in \N \cup \{ 0\} $ and $l',h,k\in \N $. 
\begin{itemize}
\item[{\em i)}] If $l\neq 0$ we have
$$f_h(l,l',k)=x_1f_{h-1}(l,l',k-1)+\cdots +x_{\min \{ l,h\} }f_{h-\min \{ l,h\} }(l,l',k-1),$$
$$f'_{h+1}(l,l',k) =  x_1f'_{h}(l,l',k-1)+\cdots 
+x_{\min \{ l,h\} }f'_{h+1-\min \{ l,h\} }(l,l',k-1)+$$ 
$$ +y_1f'_{h}(l,l',k-1)+\cdots +y_{\min \{ l',h\} }f'_{h+1-\min \{ l',h\} }(l,l',k-1)+$$ $$ 
+y_1f_{h}(l,l',k-1)+\cdots +y_{\min \{ l',h+1\} }f_{h+1-\min \{ l',h+1\} }(l,l',k-1)
.$$
\item[{\em ii)}] Let $k\neq 1$. If $h\leq k-1$ then $f_h(l,l',k)$, $\  f'_h(l,l',k)=0$. If 
$h>k-1$\linebreak and $l\neq 0$ then $f_h(l,l',k)$ is a polynomial of $\Z 
[x_1,\ldots ,x_{\min \{ l,h-k+1\} }]$ \linebreak homogeneous of degree $k$ 
and $f'_h(l,l',k)$ is a polynomial of \linebreak $\Z [x_1,\ldots ,x_{\min \{ 
l,h-k+1\} },y_1,\ldots ,y_{\min \{ l',h-k+1\} }]$ homogeneous of 
degree $k$.
\item[{\em iii)}] If $p,\  m'\in \N $, $\  p$ is prime, $p>k$ and $m'$ is a coefficient
of either one of $f_h(l,l',k)$, $\  f'_h(l,l',k)$ or of the polynomial 
obtained by setting $y_1=1$ 
and $y_2,\ldots ,y_{l'}=0$ in $f'_h(l,l',k)$ then 
$p\; \nmid \; m'$.
\end{itemize}\end{lemma}
\pf i) follows from the definitions. If 
$k\neq 1$,
$\  j\in \N \cup 
\{ 0\} $ and $h-(k-1)<j\leq h$ we have $f_0(l,l',k-1)$, $f_{h-j}(l,l',k-1)$, $f'_{h-j}(l,l',k-1)=0$, hence 
by i) and induction on $k$ we get ii). If 
$l=0$ then iii) may be proved by induction on $l'$ using 
the relation
$$F(0,l',k)=\sum _{j=0}^k{k\choose
j}(y_{l'})^{k-j}y^{l'(k-j)}F(0,l'-1,j).$$
If $l\neq 0$ then iii)  may be proved by induction on $l$ using 
the relation
$$F(l,l',k) = \sum _{j=0}^k{k\choose
j}(x_l)^{k-j}x^{l(k-j)}F(l-1,l',j)(x,0)+$$ $$ +
\sum _{j'=1}^k{k\choose
j'}(x_l)^{k-j'}y^{l(k-j')}F(l-1,l',j')(0,y).\  \square  
$$
\newline \newline 
For $l=2,\ldots
,\m $ let $\Phi (n,l)\in \cM ((n-l+2)\times n, \Z [x_1,\ldots , 
x_{n-l+1}])$
be such that the entry of $\Phi (n,l)$ of indices $h,k$ 
is $f_{h-1}(n-l+1,l',n-k)$ 
for  $h=1,\ldots ,n-l+2$, $\  k=1,\ldots
,n$, where $l'$ is any element of $\N $.\newline
For $l=2,\ldots
,\m $ let $\  \Psi (n,l)\in \cM ((l-1)\times n,\Z [x_1,\ldots ,x_{l-2},
y_1,\ldots ,y_{l-1}])$  
be such that the entry of  $\Psi (n,l)$ of indices  $h',k$ is 
$f'_{h'}(l',l-1,n-k)$  for $ h'=1,\ldots ,l-1$, $\  k=1,\ldots ,n$, where
$ l'$ is any element of $\N 
\cup \{ 0\} $ such that $l'\geq l-2$.  
\vspace{4mm}\newline 
{\bf Examples.}
We
have that $\Phi (7,3)\in \cM (6\times 7,\Z [x_1,\ldots ,x_5])$ is the following matrix:
\vspace{2mm}
$$\pmatrix{0 &
  0 & 0 & 0 & 0 & 0 & 1\cr
  & & & & & & \cr
  0 & 0 & 0 & 0 & 0 & x_1 & 0 \cr
 & & & & & & \cr
 0 & 0 & 0 & 0 & (x_1)^2 & x_2 & 0 \cr
 & & & & & & \cr 0 & 0 & 0 & (x_1)^3 & 2x_1x_2 & x_3 & 0 \cr
  & & & & & & \cr
  0 & 0 & (x_1)^4 &
3(x_1)^2x_2 &
(x_2)^2+2x_1x_3 & x_4 & 0 \cr
   & & & & & & \cr
 0 & (x_1)^5 & 4(x_1)^3x_2 & 3(x_1)^2x_3+3x_1(x_2)^2 &
2x_1x_4+2x_2x_3 & x_5 & 0 \cr }
\vspace{3mm}$$ 
and  $\Psi (7,4)\in \cM (3\times 7, \Z [x_1,x_2,y_1,y_2,y_3])$ is the following matrix:\vspace{2mm}  
$$\pmatrix{ 0 & 0 & 0 & 0 & 0 & y_1 & 0 \cr &  &  &  &
 & & 
\cr 0 & 0 & 0 & 0 & (y_1)^2 +2x_1y_1 & y_2 
& 0 \cr   & & & & & &  \cr 0 & 0 & 0 &
\begin{array}{c}(y_1)^3+\\+3x_1(y_1)^2+3(x_1)^2y_1\end{array} & 
\begin{array}{c}2y_1y_2+\\ +2x_1y_2+2x_2y_1\end{array} & y_3 & 0  \cr
}.\vspace{4mm}$$
For $i=0,\ldots ,n-1$ let $\Phi (n,l,i)\in \cM 
((n-l+2)\times (n-i),\Z [x_1,\ldots ,x_{n-l+1-i}])$ 
be defined as follows: if $i>n-l+1$ it is the zero matrix, otherwise it is 
the matrix obtained by putting $i$ rows more of zeros on the top 
of the submatrix of $\Phi (n,l)$ obtained by taking the first $n-l+2-i$
rows and the last $n-i$ columns. For $i=0,\ldots ,n-1$ let  
$\Psi (n,l,i)\in \cM ((l-1)\times (n-i),\Z [x_1,\ldots ,x_{l-2-i},y_1,\ldots ,y_{l-1-i}])$
be defined as follows: if $i>l-2$ it is the zero matrix, otherwise it is 
the matrix obtained by putting $i$ rows more of
zeros on the top of the submatrix of $\Psi (n,l)$ obtained by taking 
the first $l-1-i$ rows and the last $n-i$ columns.
Let $N={\displaystyle \frac{n(n+1)}{2}}$ and let us define the following 
matrix of $\cM ((n+1)\times N,\Z 
[x_1,\ldots ,x_{n-l+1},y_1,\ldots ,y_{l-1}])$:
$$\Upsilon '(n,l)=\pmatrix{\Phi (n,l,n-1) & \cdots & \Phi (n,l,0)\cr  &  &  \cr 
\Psi (n,l,n-1) & \cdots & 
\Psi (n,l,0)\cr };$$
let $\Upsilon ''(n,l)$ be the submatrix of $\Upsilon '(n,l)$ 
obtained by taking the first $n$ rows and let 
$\Upsilon (n,l)\in \cM (n\times N,\Z [x_1,\ldots 
,x_{n-l+1}, y_1,\ldots ,y_{l-1},z])$ be the matrix obtained
by adding to the row of index $n-l+2$ of $\Upsilon ''(n,l)$ the last row 
of $\Upsilon '(n,l)$ 
multiplied by $z$.\newline
Let
$$L_n=\{ (i,j)\in \{ 0,\ldots ,n-1\} \times \{ 0,\ldots ,n-1\} \ 
:\ i+j<n\} ,$$
which has cardinality $N$, and let us consider 
in $L_n$ the lexicographic order. We will consider any subset of $L_n$ 
as an ordered subset.  
\begin{lemma}\label{m}
\begin{itemize}
\item[{\em 1)}] If $A\in N(n,K)$, $\  \alpha _1,\ldots ,\alpha _{n-1}\in K$,
$\ B=\alpha _1A+\cdots
+\alpha _{n-1}A^{n-1}$ and $\ v\in K^n$ is such that
$A^{n-1}v\neq 0$ then
for $(i,j)\in L_n$
the entries of the column of $\Upsilon (n,2)(\alpha _1,\ldots ,\alpha _{n-2},0,1,\alpha _{n-1})$
of indices $(i,j)$ are the coordinates of $A^iB^jv$ with respect to the basis 
$\{ v,
Av,\ldots ,A^{n-1}v\} $.
\item[{\em 2)}] If $l\in \Big\{ 3,\ldots ,\m \Big\} $, $\  \alpha _1, \ldots ,
\alpha _{n-l+1},\beta _1, \ldots ,\beta _{l-2}\in K $, $\ c\in K\setminus \{ 0\} $,
$$B=\alpha _1A_{l,n}(c)+\cdots +\alpha _{n-l+1}(A_{l,n}(c))^{n-l+1}+\beta _1J_{l,n}+\cdots 
+\beta _{l-2}
(J_{l,n})^{l-2}$$ 
and $v\in K^n$ is such that $(J_{l,n})^{l-1}v\neq 0$ 
then for
$(i,j)\in L_n$ the entries of the column of
$\Upsilon (n,l)(\alpha _1,\ldots ,\alpha _{n-l+1},
\beta _1,\ldots ,\beta _{l-2},0,c^{-1})$ of indices $(i,j)$ are the coordinates of
$(A_{l,n}(c))^iB^jv$ with respect to the basis
$$\{ v,A_{l,n}(c)v,\ldots ,(A_{l,n}(c))^{n-l+1}v,J_{l,n}v,\ldots ,(J_{l,n})^{l-2}v\} .$$
\end{itemize}\end{lemma}
\pf By the expression of $B$ as linear combination of powers of
$A$ we get 1).  We have $A_{l,n}(c)J_{l,n}=(J_{l,n})^2$
and $(J_{l,n})^{l-1}=c^{-1}(A_{l,n}(c))^{n-l+1}$. Hence by the expression of
$B$ as linear combination of powers of $A_{l,n}(c)$ and $J_{l,n}$ we get
2). $\square $
\newline \newline 
Let $\Pi _n(L_n)$ be the set of all the subsets of $L_n$ 
of cardinality $n$. If $\Xi $ is an $n\times N$ matrix let $L_n$ 
be the set of the indices of the columns of $\Xi $ and for  
$I\in \Pi _n(L_n)$ let $\mu (I,\Xi )$ be the  
minor of $\Xi $ having $I$ as set of the indices of the columns. 
\begin{proposition}\label{4.3}
If $n\in \N $, $m\in \Big\{ 3,\ldots ,\m \Big\} $ and either $\char K=0$ or 
$\char K\geq \n $ the variety
$\widetilde {\cH }_m(n,K)$ is irreducible.
\end{proposition} \pf If $l\in \N $ and $\n +1<l\leq n-m+3$ the 
subset of $\widetilde {\cH }_m(n,K)$ of all $(A,B)$ such that $\ind A=l$ is 
irreducible, since  if $A\in \cN (n,K)$ the 
variety $\cN _A$ is irreducible. Hence, by induction on $\m -m$, it is sufficient to prove 
that if
$(\overline A,\overline B)\in \widetilde 
{\cH }_m(n,K)$, $\  \ind \overline A=n-m+2$ and $\cA $ is an open 
subset of $\widetilde {\cH }_m(n,K)$ such that $(\overline A,\overline 
B)\in \cA $ then
$$\cA \;  \cap \;  \{ (A,B)\in \widetilde {\cH }_m(n,K)\ : \    
\ind  A=n-m+3\}  \neq \emptyset . $$
By Lemma \ref{4.1} we may assume that there exists $\overline f\in K[x]$ 
such that $x\;  \mid \;  \overline f$ and  $\  \overline B
+\overline f(\overline A)$ is conjugate to $J_{m,n}$.
Since the map from $\widetilde {\cH }_m(n,K)$ into itself defined 
by $(A,B)\mapsto (A,B+\overline f(A))$ is an 
automorphism of $\widetilde {\cH }_m(n,K)$, 
we may assume that $\overline B$ is conjugate to $J_{m,n}$. Hence 
by a) of Lemma \ref{4.2} we 
may assume that 
there exists $\overline c\in K\setminus 
\{ 0\} $ such that $\overline A=A_{m,n}(\overline c)$ and there exist 
$\overline {\gamma }_1,\ldots ,\overline {\gamma }_{m-1}\in K$ 
such that $\overline {\gamma }_1\neq 0$ and
$\overline B=\overline {\gamma }_1J_{m,n}+ \overline{\gamma 
}_2\overline
A\, \overline B+
\cdots +\overline {\gamma }_{m-1}(\overline A) ^{m-2}\overline B$.\newline
In the ring $K[x,y]/(x^{n-1}y)$ let $\xi 
=(x^{n-1}y)+x$, $\  \eta =(x^{n-1}y)+y$, $\ \eta '=\eta -\overline 
{\gamma }_2\xi \eta -\cdots -\overline {\gamma }_{m-1}\xi 
^{m-2}\eta $. We have
$$\xi ^i \eta =\xi ^i \eta '+\overline {\gamma }_2\xi ^{i+1}\eta 
+\cdots +\overline {\gamma }_{\min \{ m-1,n-1-i\} }\xi ^{\min \{ 
m-2+i,n-2\} }\eta $$
for $i=1,\ldots ,n-2$, hence $\xi ^{n-2}\eta =\xi ^{n-2}\eta '$ 
and by reverse induction on $i$ we may get an expression of $\xi ^i\eta $ as 
linear 
combination over $K$ of $\xi ^i \eta ', \ldots ,\xi ^{n-2}\eta '$ for 
$i=1,\ldots ,n-2$. 
Then the map from $\widetilde {\cH }_m(n,K)$ into itself defined by 
$$(A,B)\mapsto (A,(\overline {\gamma }_1)^{-1}(B-\overline 
{\gamma }_2AB-\cdots -\overline {\gamma }_{m-1}A ^{m-2}B))$$ 
is also an automorphism of $\widetilde {\cH }_m(n,K)$. Hence we may 
assume $\overline B=J_{m,n}$.  \newline
Let  
$$\widehat {\cH }_m(n,K)=\{ (A,B,v)\in \widetilde {\cH }_m (n,K)\times K ^n \  :\  \dim
\langle  A^iB^jv\  :\  (i,j)
\in Z_m\rangle =n \} $$
and let $\widehat {\cA } =\{ (A,B,v)\in \widehat {\cH }_m(n,K)\ :\ (A,B)\in 
\cA \} $. Let $\overline v\in K^n$ be such that $(J_{m,n})^{m-1}\overline v\neq
0$. 
By b) of Lemma \ref{4.2} we have $(A_{m,n}(\overline c)
,J_{m,n},\overline v)\in \widehat {\cA }$. 
It is sufficient to prove that 
$$\widehat {\cA }\;  \cap \;  \{ (A,B,v)\in \widehat {\cH }_m(n,K)\  : \ 
\ind  A=n-m+3\} \neq \emptyset .$$  
\newline
The group GL$(n,K)$ acts on $\widehat {\cH }_m(n,K)$ by the relation
$$G\cdot (A,B,v)=(G^{-1}AG,G^{-1}BG,G^{-1}v)$$
for $(A,B,v)\in \widehat {\cH
}_m(n,K)$ and $G\in $ GL$(n,K)$.
We may assume that $\widehat {\cA }$ is stable
with respect to this action.\newline 
Let $\cG (N-n,K)$ be the grassmannian of all the
subspaces of
dimension $N-n$ of the $N$-dimensional vector space 
$K[x,y]/(x,y)^n$ 
over $K$. 
Let $\cF _m(N-n,K)$ be the open subset of $\cG (N-n,K)$ of all  
$V$ such that $V\cap \langle (x,y)^n+x^iy^j\ :\ (i,j)\in Z_m\rangle =\{ 0\} $
and let $\cI _m(N-n,K)$ be the subset of $\cF _m(N-n,K)$
of all $V$ which are ideals of the ring $ K[x,y]/(x,y)^n$.\newline
Let
$$\zeta _m\ : \  \widehat {\cH }_m(n,K)\to \cG (N-n,K)$$ 
be the map which associates to $(A,B,v)\in \widehat {\cH }
_m(n,K)$ the kernel of the homomorphism 
from $K[x,y]/(x,y)^n$ onto $K^n$ defined by $(x,y)^n+g\mapsto g(A,B)v$
for $g\in K[x,y]$.  
If $(A,B,v)\in \widehat {\cH }_m(n,K)$ 
all the elements of the 
orbit of $(A,B,v)$ have the same image under $\zeta _m$
and
the matrix having as columns the vectors
$A^iB^jv$ for $(i,j)\in L_n$ is the matrix of the coefficients of
a system of $n$ linear homogeneous equations for $\zeta 
_m((A,B,v))$. Hence $\zeta _m$ is a morphism. Moreover $\zeta _m(\widehat {\cH }_m(n,K))
\subseteq \cI _m(N-n,K)$. 
\newline
For $V\in \cF _m(N-n,K) $ let $\alpha _V$ and $\beta _V$ be the endomorphisms
of the vector space $K[x,y]/V$ over $K$ such that 
$\alpha _V(V+x^iy^j)=V+x^{i+1}y^j$ 
and $\beta
_V(V+x^iy^j)=V+x^iy^{j+1}$ for $(i,j)\in Z_m$. 
Let $A_V$ and $B_V$ be 
the matrices of $\alpha _V$ and $\beta _V$ respectively with respect to the
basis $\{ V+x^iy^j\  : \ (i,j)\in Z_m\} $, which has the vector $V+1$ as $n$-th vector.
Let $$\widehat {\zeta }_m\  :\ 
\cF
_m(N-n,K)\to M(n,K)\times M(n,K)\times K^n$$ be defined by $\widehat 
{\zeta }
_m(V)=(A_V,B_V,e_n)$ for $V\in \cF _m(N-n,K)$.\newline 
Let $V\in \cF _m(N-n,K)$ and let $\Theta $ be the
matrix of the coefficients of a system of $n$ linear homogeneous
equations for $V$. Let $(i',j')\in L_n $
and let $\Theta ^{(i',j')}$ be the column of $\Theta $ of index $(i',j')$. 
Then the $n$-tuple of the coordinates of
 $V+x^{i'}y^{j'}$
with respect to the basis $\{ V+x^iy^j\  :\  (i,j)\in Z_m\} $ is
the solution of the Cramer system
$$\mu (Z_m,\Theta )X=\Theta ^{(i',j')}$$
where $X$ is the column of the unknowns. 
Hence $\widehat {\zeta }_m $ is a morphism.\newline 
If $V\in \cI _m(N-n,K)$ we have $\alpha _V(V+g)=V+gx$ and $\beta 
_V(V+g)=V+gy$ for $g\in K[x,y]$, hence $(A_V,B_V,e_n)\in \widehat 
{\cH }_m(n,K)$ and $\zeta _m((A_V,B_V,e_n))=V$. If $(A,B,v)\in 
\widehat {\cH }_m(n,K)$ and $V=\zeta _m((A,B,v))$ we have that
$A_V$, $B_V$ are the matrices which represent $A,\ B$ respectively 
with respect to the basis $\{ A^iB^jv\ :\ (i,j)\in Z_m\} $, hence $(A,B,v)$ and $
(A_V,B_V,e_n)$ belong to the same orbit.  Then, since
$$\zeta _m(\widehat {\cA })=\big( \widehat {\zeta }_m\mid _{\cI _m 
(N-n,K)}\big) ^{-1}(\widehat {\cA }),$$
we have that $\zeta _m(\widehat {\cA })$   
is an open subset of  $\cI _m (N-n,K)$.
\newline 
Let $\widehat {\cH }'_m(n,K)$ be the subset of $\widehat {\cH }_m(n,K)$ 
of all $(A,B,v)$ such that either $m=3$ and $A$ is regular or
$m>3$ and 
there exists $G\in $ GL$(n,K)$, $\ c\in K\setminus \{ 0\} $ and $\alpha _1,\ldots ,
\alpha _{n-2m+4}\in K$ such that $G^{-1}AG=A_{m-1,n}(c)$ and
$$G^{-1}BG=\alpha _1A_{m-1,n}(c)+\cdots +\alpha _{n-2m+4}(A_{m-1,n}(c))^{n-2m+4}+J_{m-1,n}.$$
We want to prove that
$$\zeta _m(\widehat {\cA })\; \cap \; \zeta _m(\widehat {\cH 
}'_m(n,K))\neq \emptyset .$$
Let $\overline {\Upsilon }(n,m)$ be the matrix over $K$ obtained by 
setting $x_1,\ldots ,x_{n-m+1},\linebreak y_2,\ldots ,y_{m-1}=0$, $y_1=1$ and 
$z=\overline c^{-1}$ in $\Upsilon (n,m)$. Let $\widehat {\Upsilon }(n,m-1)$ be 
the matrix over $\Z [x_1,\ldots ,x_{n-2m+4},z]$ obtained by setting 
$x_{n-2m+5},\ldots ,x_{n-m+2},\linebreak y_2,\ldots ,y_{m-2}=0$ 
and $y_1=1$  in $\Upsilon 
(n,m-1)$. 
By using the definition of $\Upsilon (n,m)$, $\ \Upsilon (n,m-1)$ 
and Lemma \ref{A} we get the following properties for the 
determinants of the minors of order $n$ of $\overline {\Upsilon 
}(n,m)$ and $\widehat {\Upsilon }(n,m-1)$. \newline
We have $$\det \mu (Z_{m-1},\widehat {\Upsilon }(n,m-1))=1, \quad \det \mu 
(Z_m,\widehat {\Upsilon }(n,m-1))=z.$$ For  $i\in \{ 1,\ldots 
,n-2m+4\} $ let $Z_{m-1,i}$ be the set obtained by substituting 
$(n-m+2,0)$ with $(n-m+2-i,1)$ in $Z_{m-1}$; then we have $$\det \mu (Z_{m-1,i}, 
\widehat {\Upsilon }(n,m-1))=x_{i}.$$ \newline
Let $I\in \Pi _n(L_n)$. 
If $l\in \Big\{ 2,\ldots ,\m \Big\} $ we will say that $I$  has the property 
$\cP (l)$ if 
there exist $h_i\in \{ 0,\ldots ,i-1\} $ for $i=1,\ldots ,l-2$ such that 
$$I= \{ (n-l,0),\ldots ,(0,0)\} \; \cup \; \{ (h_i,i-h_i)\ :\   i=1,\ldots ,l-2\} 
\;  \cup \;   \{ (h,k)\} $$
where either $(h,k)=(n-l+1,0)$ or 
$h\in \{ 0,\ldots ,l-2\} $, $k=l-1-h$. \newline
If $I$ hasn't the property $\cP (m-1)$ then $\det \mu (I, \widehat 
{\Upsilon }(n,m-1))$ hasn't a term different from 0 and of degree 
0 with respect to $x_1,\ldots ,x_{n-2m+4}$. If $I$ has the 
property $\cP (m-1)$ there exist 
$j_I\in \{ 0,1\} $, $\  f_I\in \Z 
[x_1]$ and $F_I\in \Z [x_1,\ldots ,x_{n-2m+4}]$ such that $x_1\; \mid 
\; f_I$ and
$$\det \mu (I,\widehat {\Upsilon }(n,m-1))=z^{j_I}(1+f_I)+F_I;$$
if $I\neq Z_m,Z_{m-1}$ and either $\char K=0$ or $\char K>m-1$ by ii) and iii) 
of Lemma \ref{A} we 
have $f_I\neq 0$.\newline
We have $\det \mu (I,\overline {\Upsilon }(n,m))\neq 
0$ if and only if $I$ has the property $\cP 
(m)$. Hence $\det \mu (Z_{m-1},\overline {\Upsilon }(n,m))
=0$ and $\det \mu (Z_{m-1,i},\overline {\Upsilon }(n,m))=0$ for $i=1,\ldots 
,n-2m+4$.
\newline
If $I$ has the property $\cP (m)$ and $(n-m+1,0)\in I$ then $I$ 
has the property $\cP (m-1)$ and $j_I=1$. \newline
Let us assume that $I$ has the property $\cP (m)$ and there exists 
$h\in \{ 0,\ldots ,m-2\} $ such that $(h,m-1-h)\in I$. Let 
$h_{m-2}\in \{ 0,\ldots ,m-3\} $ be such that 
$(h_{m-2},m-2-h_{m-2})\in I$. Let $\widehat M_I$ be the 
minor of order 2 of $\widehat {\Upsilon }(n,m-1)$ 
obtained by taking the rows of indices $n-m+2,\  n-m+3$ 
and the columns of indices $(h,m-1-h),\  (h_{m-2},m-2-h_{m-2})$. 
Then there exists $g_I\in \Z [x_1]$ such that $x_1\;  \mid \;  g_I$ and
$$\det \mu (I,\widehat {\Upsilon }(n,m-1))=(1+g_I)\det \widehat M_I.$$
Moreover  there exists $h_I\in \Z [x_1]$ and $H_I\in  \Z [x_1,\ldots 
,x_{n-2m+4}]$ such that $x_1\; \mid \; h_I$ and the entry of $\widehat 
M_I$ of indices $(2,2)$ is $z(1+h_I)+H_I$. 
The other entries of $\widehat M_I$ are polynomials of $\Z 
[x_1,\ldots ,x_{n-2m+4}]$.
Since for $h\in \{ 0,\ldots ,m-2\} $ we have 
$(n-m+1-h)-(m-1-h)+1=n-2m+3$, by ii) and iii) of Lemma \ref{A}
if either 
$\char K\neq 0$ or $\char K>m-1$  the 
entry of $\widehat M_I$ of indices $(1,1)$ is an homogeneous polynomial of 
$\Z [x_1,\ldots ,x_{n-2m+3}]$ of degree $m-1-h$ if $h\neq m-2$ and is $x_{n-2m+3}$ if
$h=m-2$; hence $\det \mu 
(I,\widehat {\Upsilon }(n,m-1))$ has degree greater than 1. 
\newline
By Lemma \ref{m} we have that $\overline {\Upsilon 
}(n,m)$ is the matrix of the coefficients of a system of $n$ 
linear homogeneous equations for $\zeta _m((A_{m,n}(\overline 
c),J_{m,n},\overline v))$ and if $c\in K\setminus \{ 0\} $, 
$\  \alpha _1,\ldots ,\alpha _{n-2m+4}\in K$ and $v\in K^n$ is such that 
$(J_{m-1,n})^{m-2}v\neq 0$ then 
 $\widehat 
{\Upsilon }(n,m-1)(\alpha_1,\ldots ,\alpha _{n-2m+4},c^{-1})$ is 
the matrix of the coefficients a system of $n$ linear homogeneous 
equations for
$\zeta _m((A_{m-1,n}(c),\alpha _1A_{m-1,n}(c)+\cdots +\alpha 
_{n-2m+4}(A_{m-1,n}(c))^{n-2m+4}+J_{m-1,n},v))$.
\newline
Then if either $\char K=0$ or $\char K>m-1$ we get that the coordinates of 
$\zeta _m((A_{n,m}(\overline c),J_{m,n},\overline v))$ as element of 
$\cG (N-n,K)$ satisfy the equations of 
the closure in $\cG (N-n,K)$ of $\zeta _m(\widehat {\cH }'_m(n,K))$, and hence 
$\widehat {\cA }\  \cap \  \widehat {\cH }'_m (n,K) \neq \emptyset $.
This proves the result.
$\square $
\newline  \newline
{\bf Remark.} There exist $n,\ m,\ p\in \N $ such that $\ m\in
\Big\{ 3, \ldots , \m \Big\} $, $p$ is prime, $p\leq m-1 $ and if $\char K=p$ 
the variety $\widetilde {\cH
}_m(n,K)$ is reducible (and hence the variety $\cH _m(n,K)$ is reducible). 
In fact, if $\char K=3$ and $$I=\{ (3,0),(2,0),(1,0),(0,3),(0,2),(0,1),(0,0)\} $$
we have $\det \mu (I, \Upsilon (7,4))\neq 0$ but $\det \mu (I, \Upsilon (7,3))=0$; 
hence by the proof of Proposition \ref{4.3}, Lemma \ref{4.1} and Lemma 
\ref{m} if $\char K=3$ we have that $\widetilde {\cH }_4(7,K)$ is 
reducible.\newline
\begin{lemma}\label{Sh}
Let $\nu :\cY _1\to \cY _2$ be a closed surjective morphism of
quasi projective varieties and let $\cY _2$ be irreducible. If
there exists a dense subset $\widetilde {\cY }_2 $ of $\cY _2$ such that
$\overline {\nu ^{-1}(\widetilde {\cY }_2 )}=\cY _1$ and the fibers of the
elements of $\widetilde {\cY }_2 $ are all irreducible of the same dimension then
$\cY _1 $ is irreducible.
\end{lemma}
\pf It is a generalization of Theorem 8 of chap. I, \S   6 of
\cite{Sh}. $\square$ \newline  \newline
Now we can prove the result which is the aim of this paper.
\begin{theorem}\label{4.9}
If $n\in \N $ and either $\char K=0$ or $\char K\geq \n$
the variety
$\cH (n,K)$ is irreducible of dimension $n^2-1$.
\end{theorem}
\pf We prove the theorem by showing that if $(\overline
A,\overline B)\in \cH (n,K)$ and $\cA $ is an  open subset of $\cH
(n,K)$ such that $( \overline A,\overline B)\in \cA $ there exists
$(\widehat A,\widehat B)\in \cA $ such that $\widehat B$ is
regular. Let $\overline u_1,\ldots ,\overline u_{\overline t}$ be
the orders of the Jordan blocks of $\overline B$. By induction on $\overline t
$ it 
is sufficient
to show that there exists $(\widehat A,\widehat B)\in \cA $ such
that either $\widehat B$ is regular or $\rank \widehat B
>n-\overline t$.\newline
If there exists $(\widehat A,\widetilde B)\in \cA $ such
that either $\widehat A$ is regular or $\rank \widehat A
>n-\overline t$ then by the irreducibility of $\cN _{\widehat A}$ there 
exists $\widehat B\in \cN _{\widehat A }$ such that $(\widehat A, 
\widehat B)\in \cA $ 
and either 
$\widehat B$ is regular or $\rank \widehat B>n-\overline t$. 
\newline
Let $n=d_{\overline t}\overline t+ r_{\overline
t}$ where $d_{\overline t},r_{\overline t}\in \N \cup \{ 0\} $ and
$r_{\overline t}<\overline t$. Then $\ind \overline B\geq
d_{\overline t}+\min \{ 1,r_{\overline t}\} $. If $\ind \overline
B=d_{\overline t}+\min \{ 1,r_{\overline t}\} $ we have $\overline
u_1-\overline u_{\overline t}
\leq 1$, hence by Lemma \ref{4} 
and the
irreducibility of $\cN _{\overline B}$ there exists 
$\widehat A\in \cN _{\overline B}$ such that $\widehat A$ is 
regular and 
$(\widehat A,\overline
B)\in \cA $. Then, for a given
$\overline t$, we may prove the claim by induction on $\ind
\overline B$.\newline 
If there exists $i\in \{ 1, \ldots ,\overline
t-1\} $ such that $\overline u_i-\overline u_{i+1}\leq 1$ by
Proposition  \ref{4} and the irreducibility of $\cN _{\overline B}$ 
there exists $\widehat A \in \cN _{\overline B}$ such that $(\widehat A, 
\overline B)\in \cA $ and $\rank
\widehat A>n-\overline t$; hence we
may assume $\overline u_i-\overline u_{i+1}\geq 2$ for $i=1,\ldots
,\overline t-1$.\newline Let $\overline u_1-\overline u_2=2$. By
Proposition \ref{3.8'} and the irreducibility of $\cN _{\overline B}$ there
exist $\widetilde A\in \cN _{\overline B}$ and $a(\widetilde A)\in K$ such that
$(\widetilde A,\overline B)\in \cA $, $\  \rank (\overline B+a(\widetilde
A)\widetilde A) 
=n-\overline
t$, $\  \ind (\overline B+a(\widetilde A)\widetilde A)<u_1$. By the inductive
hypothesis
applied to $(\widetilde A,\overline B+a(\widetilde A)\widetilde 
A)$ and the open subset $\{ (A'',B''+a(\widetilde A)A'')\  :\  
(A'',B'')\in \cA \} $ 
there exists $(\widehat A,\widetilde B)\in \cA $ such that either $\widetilde B+a(\widetilde A)
\widehat A$ is
regular or $\rank (\widetilde B+a(\widetilde A)\widehat A)>n-\overline 
t$; hence by the
irreducibility of $\cN _{\widetilde B+a(\widetilde A)\widehat A}$ it is possible to choose 
$\widehat A$
such that either $\widehat A$ is regular or $\rank \widehat A>n-\overline 
t$. \newline 
Let $\overline
u_1-\overline u_2>2$. 
If $\overline t=2$ the claim follows by
Proposition \ref{4.3}, hence we may assume
$\overline t>2$. 
Let $\overline n=\overline u_1+ \overline
u_2$ and let $Q\in \cN (n-\overline n,K)$ be the matrix in Jordan
canonical form which has $\overline t-2$ Jordan blocks of orders
$\overline u_3,\ldots , \overline u_t$. 
Let $\cH _{\overline u_1,\overline u_3}(\overline n,K)$ be the
subvariety of $\cH (\overline n,K)$ of all $(C,D)$ such that 
$\rank D\geq \overline n-2$, $\  \,  \overline u_1\leq \ind D\leq
\overline n-(\overline u_3+1)$.
For $A,B\in \cM (n,K)$ let $A=(A_{ij})$, $\ B=(B_{ij})$ where 
$i,j=1,2$ and $A_{11}$, $\ B_{11}\in \cM (\overline n,K)$. 
Let $\cH _{\overline
u_1,\overline u_3, Q}(n,K)$ be the subvariety of $\cH (n,K)$ of
all $(A,B)$ such that $B_{12}, B_{21}=0$, $B_{22}=Q$ and
$(A_{11},B_{11})\in \cH _{\overline u_1,\overline u_3}(\overline
n,K)$. Let $$\pi _{\overline u_1,\overline u_3,Q} \ :\ \cH
_{\overline u_1,\overline u_3, Q}(n,K)\to   \cH _{\overline
u_1,\overline u_3}(\overline n,K)$$ be defined by $\pi _{\overline
u_1,\overline u_3,Q}(A,B)=(A_{11},B_{11})$. 
Since $\cH _{\overline u_1,\overline u_3,Q}(n,K)$ is homogeneous 
with respect to the entries of $A_{12}$, $A_{21}$, $A_{22}$ the 
morphism $\pi
_{\overline u_1,\overline u_3,Q} $ is closed. By Proposition
\ref{4.3} the variety  $\cH _{\overline u_1,\overline
u_3}(\overline n,K)$ is irreducible and by Lemma  \ref{2} and Lemma
\ref{3} all the fibers of $\pi _{\overline
u_1,\overline u_3,Q} $ have the same dimension. Hence by Lemma
\ref{Sh} the variety $ \cH _{\overline u_1,\overline u_3, Q}(n,K)$
is irreducible. Then there exists $(A',B')\in \cA $ such that the
Jordan blocks of $B'$ have orders $\overline n-(\overline u_3+1),
\overline u_3+1,\overline u_3,\ldots ,\overline u_{\overline t}$. If $B'$ has
this property there exists $\widehat A\in \cN _{B'}$ such that
$\rank \widehat A>n-\overline t$ and hence we get the claim. $\square $.
\newline  \newline By Theorem \ref{4.9} we get a proof of the
following result and an extension of it to algebraically closed
fields $K$ such that $\char K\geq \n $.
\begin{corollary}\label{4.10}
{\rm (J. Brian\c{c}on, 1977}). If $n\in \N $,
$\cX $ is an
algebraic surface over an algebraically closed field $K$ such that
either $ \char K=0$ or $\char K\geq \n $
and $P\in
\cX $ is nonsingular then
 $\Hilb ^n (\cO _P) $ is irreducible of dimension
$n-1$.
\end{corollary}
\pf By \cite{Na} the quotient of $\widehat {\cH } (n,K)$ with respect to the
action of \linebreak GL$(n,K)$ is a variety isomorphic to $\Hilb ^n(\cO
_P)$.
Since $\widehat {\cH }(n,K)$ is an open subset of $\cH (n,K)\times K^n$ by Theorem \ref{4.9}
it
is irreducible of dimension $n^2-1+n$. Since
the stabilizer of any element of $\widehat {\cH } (n,K)$ with respect to the
action of GL$(n,K)$ is trivial, we get the claim. $\square$

\end{document}